\documentclass[12pt]{article}
\usepackage{latexsym}
\usepackage{amsmath}
\renewcommand{\epsilon}{\varepsilon}

\usepackage{amsmath,bm}
\usepackage{amssymb}
\usepackage{ntheorem}
\usepackage[ansinew]{inputenc}
\usepackage{graphicx}
\usepackage{epsfig}
\usepackage{dsfont}
\usepackage{bbm}
\usepackage{subfig}
\usepackage{float}

\usepackage[authoryear]{natbib}
\newtheorem{satz}{Theorem}[section]

\newtheorem{rem}[satz]{Remark}
\newtheorem{assumption}[satz]{Assumption}
\newtheorem{algo}[satz]{Algorithm}

\def\3{\ss}

\def \R{I \!\! R}   %Reelle Zahlen
\def \N{I \!\! N}   %Natuerliche Zahlen
\def \Z{Z \!\!\! Z} %Gan
\def \C{\mathbb{C}}

\newcommand{\E}{\mathbbm{E}}

\newcommand{\bea}{\begin{eqnarray*}}
\newcommand{\eea}{\end{eqnarray*}}
\newcommand{\be}{\begin{eqnarray}}
\newcommand{\ee}{\end{eqnarray}}

\newcommand{\ba}{\begin{array}}
\newcommand{\ea}{\end{array}}

\def\3{\ss}

\begin{document}

\title{Testing for stationarity in multivariate locally stationary processes}

\author{ Ruprecht Puchstein,  Philip Preu\ss \\
Ruhr-Universit\"at Bochum \\
Fakult\"at f\"ur Mathematik \\
44780 Bochum \\
Germany \\
{\small email:  ruprecht.puchstein-t9p@ruhr-uni-bochum.de}\\
{\small email: philip.preuss@ruhr-uni-bochum.de}\\
}

 \maketitle
\begin{abstract}
In this paper we propose a nonparametric procedure for validating the assumption of stationarity in multivariate locally stationary time series models. We develop a bootstrap assisted test based on a Kolmogorov-Smirnov type  statistic, which tracks the deviation of the time varying spectral density from its best stationary approximation.  In contrast to all other nonparametric approaches, which have been proposed in the literature so far, the test statistic does not depend on any regularization parameters  like smoothing bandwidths or a window length, which is usually required in a segmentation of the data. We additionally show how our new procedure can be used to identify the components where non-stationarities occur and indicate possible extensions of this innovative approach. We conclude with an extensive simulation study, which shows finite sample properties of the new method and contains a comparison with existing approaches.  
\end{abstract}

AMS subject classification: 62M10, 62M15, 62G10

Keywords and phrases: goodness-of-fit test, locally stationary processes, empirical spectral measure, spectral density, integrated periodogram

\section{Introduction}
\def\theequation{1.\arabic{equation}}
\setcounter{equation}{0}
A vast amout  of scientific effort in the field of time series analysis has been conducted in the framework of second order stationarity, which assumes that the first and second order moments of the considered process are constant over time. This assumption allows for a straightforward development of statistical procedures like parameter estimation or forecasting techniques, and, for this reason,  a vast amount of academic literature exists which one can rely on in the statistical anyalsis [see for example \cite{anderson1971}, \cite{brillinger1981} or \cite{brodav1991}]. Due to this large number of existing procedures, the assumption of stationarity is usually imposed when a set of time series data is analyzed, but, since the dependency structure of real world time series very often changes over time,   this assumption frequently turns out to be too restrictive [cf. \cite{starica2005} or \cite{motivation2} among many others]. For this reason, the availability of goodness-of-fit tests, which are able to detect deviations from stationarity, is of particular importance and, in the context of parametric set-ups, there exists a vast amount of approaches to validate this condition [see \cite{Taniguchi} or \cite{changli2009} among others]. However, since these tests critically depend on the right model-choice, nonparametric procedures should be prefered and a very early approach can be found in \cite{priestleystat}.

Throughout the last decade a larger number of articles were written concerning the development of  nonparametric tests for stationarity in the framework of so called locally stationary processes. These kinds of stochastic processes were introduced by \cite{dahlhaus1997} and describe a class of non stationary processes, which can be locally approximated by stationary models in a way, which makes an asymptotic theory possible. In this setting, several  nonparametric approaches have been proposed to validate the assumption of a constant dependency structure, but essentially all of them have in common that they depend on several regularization parameters, the choice of which can heavily influence the result of the statistical analysis. Since it is in general not obvious how to choose these parameters, this is a major disadvantage from a practical point of view, which we would like to highlight in the following brief review of existing procedures.

To the best of our knowledge, the first to propose a nonparametric test for stationarity in the set-up of local stationarity were \cite{vonSachs2000}, who use a Haar wavelet series expansion of  local periodogram estimates of the spectral density function to derive a testing procedure. However, for this procedure to work well it is essential to incorporate the right amount of wavelet coefficents which is in general difficult to realize in practice. In addition, it is not obvious how the window length $N$ in the calculation of the local periodograms should be chosen. \cite{paparoditis2009} and \cite{paparoditis2010} develop goodness-of-fit tests by estimating $L_2$-distances between the local spectral density and the best approximation by a stationary spectral density  via smoothed local periodograms. While the first approach works by using quantiles of the standard normal distribution, bootstrap replications are required to obtain critical values in the second procedure. Both approaches, however, have in common that the window length $N$ and two additional smoothing bandwidths $b$ and $h$ have to be chosen.  A different idea was considered by \cite{rao2010}, who follow a Portmanteau-type approach by calculating  empirical covariances and  summing them up to a specified order $m$. However, the choice of this order is left to the practitioner and is, in general, difficult to make. Additionally, the choice of a smoothing bandwith $b$ is required in the standardization of the test statistic.  \cite{detprevet2010} estimate the same $L_2$-distance as proposed in  \cite{paparoditis2009} but avoid the additional smoothing  by calculating Riemann sums over the local periodogram estimates. This yields a procedure, which only depends on the choice of one window length and thus presents an improvement but is still far from optimal. Trying to resolve this issue, \cite{detprevet2011b} follow the proposals of \cite{dahlhaus2009} and \cite{dahlpolo2009} to measure deviations from stationarity by the Kolmogorov-Smirnov-type distance
\begin{align}
\frac{1}{2\pi}\sup_{v,\omega \in [0,1]}\Big |\int_0^{\omega\pi}\int_0^vf(u,\lambda)du d\lambda-v\int_0^{\omega\pi}\int_0^{1}f(u,\lambda)dud\lambda \Big|.\label{Distance}
\end{align}
It is easy to see that this quantity vanishes if the time-varying spectral density $f(u,\lambda)$ is independent of $u$ while it is strictly positive if $f(u,\lambda)$ depends on $u$ on a subset of $[0,1] \times [0,\pi]$ with nonzero Lebesgue measure. The authors consider two different ideas to estimate the distance measure \eqref{Distance}. The first approach works by employing a Riemann sum approximation over the local periodogram which works  well in finite sample situations but leads to the same problem as mentioned above (namely the choise of the window length). The second approach makes use of the pre-periodogram introduced by \cite{neumsach1997}. Since this approach does not require the choice of the critical window length, it solves the problem of testing for stationarity without choosing a regularization parameter in the calculation of the test statistic  from a theoretical point of view. However, the performance in finite sample situations is extremely poor as was pointed out in \cite{detprevet2011b}. The main obstacle concerning this method is that in a pre-periodogram approach all $T$ data (here and throughout the whole paper $T$ denotes the sample size) are used in the estimation of $\int_0^vf(u,\lambda)du$ regardless of how $v$ is chosen.

The above review of existing literature on the topic of testing for stationarity highlights that the derivation of a testing procedure which does not require the choice of any regularization parameters and simultaneously works  well in finite sample situations is an unsolved problem. The first goal of this paper is to resolve this issue. For this purpose, we will introduce an estimator for a quantity closely related to \eqref{Distance} which is neither based on the local periodogram nor the pre-periodogram. The main difference will be that we estimate the quantity $v^{-1} \int_0^vf(u,\lambda)du$ by the usual periodogram using the data at time-points $1,...,\lfloor vT \rfloor$. With this idea we circumvent both the segmentation issue typical for a local periodogram approach and the problem of taking into account too many data after $\lfloor vT \rfloor$ as is inherent in any pre-periodgram based procedure. Note that although this new approach is quite natural, the mathematical details are extremely difficult. This is mainly due to the fact that the discrete fourier transformation will be evaluated at frequencies corresponding to different bases which are then combined via a Riemann sum. We will show that an appropriately standardized test statistic converges to the supremum of a Gaussian process. However, since the covariance kernel of the limiting process depends on unknown quantities of the underlying data generating process, critical values for the construction of a formal test are not readily available [see \cite{dahlhaus2009} or \cite{detprevet2011b} for more details on this]. Therefore resampling methods are required to obtain critical values and for this purpose we will propose a bootstrap method that is easy to implement and works well in finite samples. It is then shown that this method can be applied to construct an asymptotic test for stationarity with a controlled type I error.

The second contribution of this paper is to treat the problem of testing for a constant dependency structure in a multivariate setting which distinguishes our method from all approaches mentioned above. In fact, although from a practical point of view it is of great interest to analyze the dependency structure between different time series (like for example stock returns or macroeconomic data), all procedures reviewed above are restricted to the investigation of the autocovariance structure of univariate time series. fTo the best of our knowledge only \cite{rao2012} consider a multivariate setting in their extension of the above mentioned Portmanteau-type test of \cite{rao2010}. However, the testing procedure developed by these authors also requires the choice of several tuning parameters and does not provide a framework for detecting the specific components, which exhibit deviations from stationarity. Methods for such a refined analysis might be of interest in practical applications, for example if one is interested in determining how many components of a data set exhibit non-stationarities and whether there exist components which are in fact stationary. We will indicate how such a refined analysis for the characterisation of non- stationarities can be conducted using  our approach. 

The remainder of this paper is organized as follows. In Section $2$ we introduce necessary notation, present the framework in which we conduct our analysis and motivate our further proceedings. In Section $3$ we define the test statistic, derive crucial asymptotic properties  and propose and validate a bootstrap algorithm to approximate its distribution under the null hypothesis. Section $4$ contains several possible extensions of our approach and a simulative study including nominal levels and estimates of the power of our procedure is presented in Section $5$. The paper concludes with an appendix in Section $6$, which contains all technical proofs.

\section{The set up}
\def\theequation{2.\arabic{equation}}
\setcounter{equation}{0}

We assume to observe a centered $\mathbb{R}^d$ valued stochastic process $\mathbb{X}(T):=\{\boldsymbol{X}_{1,T},...,\boldsymbol{X}_{T,T}\}$, with $\boldsymbol{X}_{t,T}=(X_{t,T,1},...,X_{t,T,d})^T$, where, for each $T\in\N$ and $t=1,...,T$, the random vector $\boldsymbol{X}_{t,T}$ possesses a locally stationary representation in the sense of \cite{dahlpolo2009}, i.e. there exists a MA($\infty$) representation of the form
\begin{align}\label{defprocess1}
\boldsymbol{X}_{t,T}=\sum_{l=0}^{\infty}\boldsymbol{\Psi}_{t,T,l}\boldsymbol{Z}_{t-l}, 
\end{align}
where $\{\boldsymbol{Z}_t\}_{t\in\mathbb{Z}}$ denote $d$-variate independent Gaussian random variables with some possibly time varying covariance matrix $\boldsymbol{\Sigma}_t$ and $\boldsymbol{\Psi}_{t,T,l}\in\mathbb{R}^{d\times d}$ $(l\in\mathbb{Z}$, $t=1,...,T)$ represent time varying matrices containing the linear factors of the MA($\infty$) representation. We note that the time dependence of the covariance matrices $\boldsymbol{\Sigma}_t$ of the innovation process can be  included in the linear factors by replacing the matrices  $\boldsymbol{\Psi}_{t,T,l}$ by $\boldsymbol{\Sigma}_{t-l}^{1/2}\boldsymbol{\Psi}_{t,T,l}$, and it is therefore not restrictive to assume  that the random vectors $\boldsymbol{Z}_t$ have unit covariance matrix.

The following set of assumptions ensures that $\boldsymbol{X}_{t,T}$ can be locally approximated by a stationary process and is rather standard in the framework of local stationarity; cf. \cite{dahlhaus2000}, \cite{dahlpolo2009} or \cite{paparoditis2009} among many others.

\begin{assumption}
\label{ass1}
We assume that for each $T\in\N$ and each $t\in\{1,...,T\}$ the random variable $\boldsymbol{X}_{t,T}$ has a MA($\infty$) representation  \eqref{defprocess1}
with  independent $\boldsymbol{Z}_t\sim N(0,\boldsymbol{I}_d)$ such that there exist twice continuously differentiable functions $\boldsymbol{\Psi}_l:[0,1]\rightarrow\mathbb{R}^{d\times d}$ $(l\in \N)$ which satisfy the following set of conditions:
\begin{itemize}
\item[i)] The linear factors $\boldsymbol{\Psi}_{t,T,l}$ can be approximated by the functions $\boldsymbol{\Psi}_l$ such that
\be \label{apprbed}
\sum_{l=0}^{\infty}\sup\limits_{t=1,...,T}\|\boldsymbol{\Psi}_l(\frac{t}{T})-\boldsymbol{\Psi}_{t,T,l}\|_{\infty}=O(\frac{1}{T}).
\ee
\item[ii)] The functions $\boldsymbol{\Psi}_l$ fulfill  the conditions
\be
\sum_{l=0}^{\infty}&\sup\limits_{u\in[0,1]}\|\boldsymbol{\Psi}_l(u)\|_{\infty}|l|<\infty,\label{summ1}\\
\sum_{l=0}^{\infty}&\sup\limits_{u\in[0,1]}\|\boldsymbol{\Psi}_l'(u)\|_{\infty}|l|<\infty,\label{summ2}\\
\sum_{l=0}^{\infty}&\sup\limits_{u\in[0,1]}\|\boldsymbol{\Psi}_l''(u)\|_{\infty}<\infty,
\ee
where $\| \cdot \|_\infty$ denotes the usual maximum norm of a matrix.
\end{itemize}
\end{assumption}

It should be noted that the assumption of Gaussianity can be abandoned and that the extension of the results to non-Gaussian innovations is cumbersome but straigthforward; see Remark \ref{gaussian} for more details. It is furthermore straightforward to show that a large class multivariate tvARMA(p,q) processes is included in this theoretical framework; cf. the proof of Proposition 2.3 in \cite{dahlpolo2009} for more details on this. For a process $\mathbb{X}(T)=\{\boldsymbol{X}_{t,T}\}_{t=1,...,T}$ with  representation of the form \eqref{defprocess1}, one can define the $\C^{d \times d}$ valued time-varying spectral density matrix by
\begin{equation}
\boldsymbol{f}(u,\lambda):=\frac{1}{2\pi}\sum_{l,m=0}^\infty\boldsymbol{\Psi}_{l}(u)\boldsymbol{\Psi}_{m}(u)^T \exp(-i\lambda(l-m))
\end{equation}
[see \cite{dahlhaus2000}], which does not depend on $u$ if the underlying process is stationary. Therefore, in order to develop a nonparametric test for stationarity, it is natural to investigate whether $\boldsymbol{f}(u,\lambda)$ is constant in $u$-direction  [see also \cite{vonSachs2000}, \cite{paparoditis2009} or \cite{paparoditis2010} among many others]. For this reason we define a Kolmogorov-Smirnov type distance by

\bea
D&:=&\|[\sup\limits_{(v,\omega)\in[0,1]^2}|[\boldsymbol{D}(v,\omega)]_{a,b}|]_{a,b=1,...,d} \|_F
\eea
where $\| \cdot \|_{F}$ denotes the Frobenius norm and the matrix $\boldsymbol{D}(v,\omega)\in\C^{d \times d}$ is defined by
\bea
\boldsymbol{D}(v,\omega)&:=&\frac{v}{2\pi}\Big(\int_0^{\omega\pi}\int_0^v\boldsymbol{f}(u,\lambda)du d\lambda-v\int_0^{\omega\pi}\int_0^{1}\boldsymbol{f}(u,\lambda)dud\lambda\Big).
\eea
If $\mathbb{X}(T)$ is stationary, then $D$ obviously vanishes, while it is strictly positive if  $\boldsymbol{f}(u,\lambda)$ is non-constant in $u$-direction on a  set with positive Lebesgue measure. Thus, for developing an asymptotic level $\alpha$-test for the null hypothesis
\be \label{null}
\text{H}_0: \quad \boldsymbol{f}(u,\lambda)\text{ is independent of }u
\ee
versus
\be \label{H1}
\text{H}_1: \quad \boldsymbol{f}(u,\lambda)\text{ is not constant in }u \text{-direction on a set with positive Lebesgue measure},
\ee
we have to
\begin{itemize}
\item[(1)] define an appropriate estimator $\hat D_T$ for $D$,
\item[(2)] obtain an estimator $\hat q_{1-\alpha}$ for the $(1-\alpha)$-quantile $q_{1-\alpha}$ of $\hat D_T$ under $H_0$,
\item[(3)] reject $H_0$ if $\hat D_T$ surpasses $\hat q_{1-\alpha}$.
\end{itemize}
Note that, in the one-dimensional case, the quantity $D$ is different from \eqref{Distance} since the entries of the matrices $\boldsymbol{D}(v,\omega)$ are weighted down to zero if $v$ approaches the origin, which is necessary, because the estimate of $\int_0^vf(u,\lambda) du$, which we will present in the next section, is based on  the first $2\lfloor vT/2 \rfloor$ data, resulting in an increasing variance for $v \rightarrow 0$.

We now proceed as follows: We  define an estimator $\hat D_T$ for the distance measure $D$ in Section \ref{Teststatsection} and propose a valid algorithm for estimating the critical values $q_{1-\alpha}$ in Section \ref{Bootstarpsection}. We will then show that by following the above three steps, we obtain a consistent asymptotic level-$\alpha$ test for the null hypothesis of a constant dependency structure.

\section{Testing for stationarity}
\def\theequation{3.\arabic{equation}}
\setcounter{equation}{0}
\subsection{The test statistic $\hat D_T$}
\label{Teststatsection}
In order to obtain an empirical version $\hat D_T$ of $D$ we have to define an estimator for $\int_0^{v}\boldsymbol{f}(u,\lambda)du$,  $v\in[0,1]$, at first. For this reason, we take  $2\lfloor vT/2\rfloor$, which is the even integer closest to $vT$, and calculate the usual periodogram  based on the first $2\lfloor vT/2\rfloor$ data points, i.e. we compute
\begin{align*}
 \boldsymbol{I}_{2\lfloor vT/2\rfloor}(\lambda):=\frac{1}{4 \pi \lfloor vT/2\rfloor} \sum_{r,s=0}^{2\lfloor vT/2\rfloor-1} \boldsymbol{X}_{1+s,T}\boldsymbol{X}_{1+r,T}^T \exp(-i \lambda (s-r)).
\end{align*}
The proof of Theorem \ref{hauptsatz} will reveal that this results in an   asymptotically unbiased estimator for $v^{-1} \int_0^{v}f(u,\lambda)du$, thus an estimator for $D$ and $\boldsymbol{D}(v,\omega)$ is  obtained by calculating

\be
\label{supstat}
\hat D_T&:=&\|[\sup\limits_{(v,\omega)\in[0,1]^2}|[\hat {\boldsymbol D}_T(v,\omega)]_{a,b}|]_{a,b=1,...,d} \|_F
\\
\label{empprocess}
\hat {\boldsymbol D}_T(v,\omega)&:=&v\Big(\frac{1}{T}\sum_{k=1}^{\lfloor \omega   \lfloor vT/2\rfloor   \rfloor}\boldsymbol I_{2\lfloor T/2\rfloor}\big(\lambda_{k,2\lfloor vT/2\rfloor} \big)-\frac{v}{T}\sum_{k=1}^{\lfloor \omega   T/2   \rfloor}\boldsymbol I_T\big(\lambda_{k,T}\big)\Big),
\ee

where $\lambda_{k,n}=2\pi k/n$  $(n\in\N ,k=1,...,n)$ denote the Fourier frequencies to the base $n$.  This means that we construct an estimator for $\boldsymbol{D}(v,\omega)$  by replacing the integral in $\lambda$ direction by a Riemann sum over the respective Fourier frequencies and substitute the estimator $2\lfloor vT/2\rfloor T^{-1}\boldsymbol{I}_{2\lfloor vT/2\rfloor}(\lambda)$ for the integral $\int_0^{v}\boldsymbol{f}(u,\lambda)du$ and $\boldsymbol{I}_{T}(\lambda)$ for the integral $\int_0^{1}\boldsymbol{f}(u,\lambda)du$ over the spectral density matrix. As mentioned above, the multiplication with the factor $v$ ensures that the estimator is weighted down when $v$ approaches zero, which is necessary, since the closer $v$ comes to the origin the less data are involved in the calculation of $\boldsymbol{I}_{2\lfloor vT/2\rfloor}(\lambda)$. The following theorem specifies the asymptotic behaviour of the process $\{\hat {\boldsymbol D}_T(v,\omega)\}_{(v,\omega)\in[0,1]^2}$ both under the null hypothesis and the alternative. Here and throughout the paper the symbol $\Rightarrow$ denotes weak convergence in $L^\infty([0,1]^2)$.
\begin{satz}\label{hauptsatz}
Suppose Assumption \ref{ass1} holds. Then we have:
\begin{itemize}
\item[a)] If \eqref{null} holds, then the process $\{\sqrt{T}\hat{\boldsymbol{D}}_T(v,\omega)\}_{(v,\omega)\in[0,1]^2}$ converges weakly to a centered Gaussian process $\{\boldsymbol{G}(v,\omega)\}_{(v,\omega)\in[0,1]^2}$, i.e.
\begin{align}\label{prozesskonvergenz}
\{\sqrt{T}\hat{\boldsymbol{D}}_T(v,\omega)\}_{(v,\omega)\in[0,1]^2}\Rightarrow\{\boldsymbol{G}(v,\omega)\}_{(v,\omega)\in[0,1]^2},
\end{align}
where the covariance kernel of $\{\boldsymbol{G}(v,\omega)\}_{(v,\omega)\in[0,1]^2}$ is given by 
\begin{align*}
\text{Cov}([\boldsymbol{G}(v_1,\omega_1)]_{a_1,b_1},[\boldsymbol{G}(v_2,\omega_2)]_{a_2,b_2})=&\frac{1}{2\pi} v_1v_2(\min(v_1,v_2)-v_1v_2)\int_0^{\min(\omega_1,\omega_2)\pi} f_{a_1b_2}(\lambda)f_{b_1a_2}(-\lambda) d\lambda
\end{align*}
for $(v_1,\omega_1),(v_2,\omega_2)\in[0,1]^2$ and $(a_1,b_1),(a_2,b_2)\in\{1,...,d\}^2$.
\item[b)] There exists a constant $C >0$, such that for all $(a,b)\in\{1,...,d\}^2$
\begin{align*}
\lim\limits_{T\rightarrow\infty}\mathbb{P}(\sup\limits_{(v,\omega)\in[0,1]^2}|[\hat D_T(v,\omega)]_{a,b}|>C)=1\quad &\text{if }\quad \sup\limits_{(v,\omega)\in[0,1]^2}|[\boldsymbol{D}(v,\omega)]_{a,b}|>0\\
\sup\limits_{(v,\omega)\in[0,1]^2}|[\hat D_T(v,\omega)]_{a,b}|=O_P(T^{-1/2})\quad &\text{if }\quad \sup\limits_{(v,\omega)\in[0,1]^2}|[\boldsymbol{D}(v,\omega)]_{a,b}|=0.
\end{align*}
\end{itemize}
\end{satz}

If one compares this new approach to construct an empirical version of ${\boldsymbol D}(v,\omega)$ with the procedure for estimating the quantity \eqref{Distance}  in \cite{detprevet2011b} it is evident that estimating the integrals $\int_0^{v}\boldsymbol{f}(u,\lambda)du$ as a whole via the local periodogram  $\boldsymbol{I}_{2\lfloor vT/2\rfloor}(\lambda)$ rather than by approximating it by means of a second Rieman sum over local periodogram estimates in $u$-direction is key to avoiding the necessity of any regularizing parameter for the blocksize. In the one dimensional case, the asymptotic variance in Theorem \ref{hauptsatz} a) is the 'same' as in Theorem 2.1 a) of \cite{detprevet2011b} with the only difference being the factor $v_1v_2$ which  results from the additional downweighting at zero. This yields that the behaviour of our new parameter free procedure  imitates the corresponding behaviour of the test statistics proposed in \cite{detprevet2011b} under local alternatives which has been investigated in more detail  by \cite{pappreuss} for three different frequency domain based tests for stationarity.

\subsection{Bootstrap procedure}
\label{Bootstarpsection}
Although  we established crucial asymptotic properties of the test statistic (\ref{supstat}) in Theorem \ref{hauptsatz}, which  characterizes its behaviour under the null hypothesis \eqref{null} and the alternative respectively, it is in general not possible to obtain the asymptotic quantiles of the test statistic $\hat D_T$ under $H_0$ from this result, since the limiting distribution in Theorem \ref{hauptsatz} depends in a complicated way on the unknown spectral density $\boldsymbol{f}$ of the data generating process [see \cite{dahlhaus2009} or \cite{detprevet2011b} for a similar situation]. In order to resolve this obstacle we now suggest an algorithmic procedure to estimate the quantiles of $\hat D_T$ under $H_0$ which is essentially an extension of the popular one dimensional AR($\infty)$ bootstrap, which was introduced by \cite{Kreiss1988}. This resampling scheme became particularly famous for being able to capture the main serial dependencies of the underlying time series quite well, and, broadly speaking, exploits the fact that a stationary process can be well approximated by an AR($p$)-model provided that the order $p$ is large enough. Before we describe the method in detail and validate its asymptotic accuracy we state the following essential assumption:
\begin{assumption}\label{annahmenbootstrap}
The stationary $\mathbb{R}^d$ valued process $\boldsymbol{X}_t$ with spectral density function $\lambda\mapsto\int_0^1 \boldsymbol{f}(u,\lambda) du$  possesses an AR($\infty$)-representation
\begin{equation}\label{AR1}
\boldsymbol{X}_t=\sum_{j=1}^{\infty}\boldsymbol{a}_j\boldsymbol{X}_{t-j}+\boldsymbol{\Sigma}^{1/2}\boldsymbol{ Z}_t,
\end{equation}
where the $\boldsymbol{Z}_t$ are independent  $N(0,\boldsymbol{I}_d)$ distributed random variables, $\boldsymbol{\Sigma} \in \R^{d \times d}$ is a positive definite matrix, and the conditions
\bea
\text{det}\Big(\boldsymbol{I}_d-\sum_{j=1}^{\infty}z^j\boldsymbol{a_j}\Big)\neq 0 \text{ for } |z|\leq 1\quad\text{ and }\quad \sum_{j=0}^{\infty}|j|\|\boldsymbol{{a}}_{j}\|_\infty<\infty
\eea
hold.
\end{assumption}
The main idea now is to choose an order $p=p(T) \in \N$ which grows to infinity for $T\rightarrow\infty$, and to fit a vector autoregressive model with order $p$ [abbr. VAR(p)] to the process \eqref{AR1}. This is done by computing a consistent estimator $(\boldsymbol{\hat{a}}_{1,p},...,\boldsymbol{\hat{a}}_{p,p})$
for the minimizer 
\begin{equation}\label{arpfitt}
(\boldsymbol{a}_{1,p},...,\boldsymbol{a}_{p,p}):=\underset{\boldsymbol{b}_{1,p},\ldots,\boldsymbol{b}_{p,p}}{\operatorname{argmin}} \text{tr} \Big(\mathbb{E}[(\boldsymbol{X}_{t}-\sum_{j=1}^p\boldsymbol{b}_{j,p}\boldsymbol{X}_{t-j})(\boldsymbol{X}_{t}-\sum_{j=1}^p\boldsymbol{b}_{j,p}\boldsymbol{X}_{t-j})^T]\Big),
\end{equation}
where the exact asymptotic requirements for this estimator will be stated in Theorem \ref{theorembootstrap}.
We then generate bootstrap replicates $\hat D_T^*$ of the test statistic $\hat D_T$ according to the following algorithmic procedure: 
\begin{algo}\label{algo1}\mbox{}
\begin{itemize}
\item[(1)] Simulate $T$ data from the model
\bea
\label{fittedarp}
\boldsymbol{X}^*_{t,T}=\sum_{j=1}^{p}\boldsymbol{\hat{a}}_{j,p}\boldsymbol{X}^*_{t-j,T}+\boldsymbol{\hat{\Sigma}}_{p}^{1/2}\boldsymbol{Z}^*_{j,T},
\eea
where $\boldsymbol{Z}^*_{j,T}$ are independent $N(0,\mathbf{I}_d)$ distributed and $\boldsymbol{\hat{\Sigma}}_{p}=\frac{1}{T-p}\sum_{j=p+1}^T(\boldsymbol{\hat{z}}_i-\boldsymbol{\bar{z}}_T)(\boldsymbol{\hat{z}}_i-\boldsymbol{\bar{z}}_T)^T$ 
with  $\boldsymbol{\bar{z}}_T:=\frac{1}{T-p}\sum_{j=p+1}^T\boldsymbol{\hat{z}}_j$ and $\boldsymbol{\hat{z}}_{j}:=\boldsymbol{X}_{j,T}-\sum_{i=1}^{p}\boldsymbol{\hat{a}}_{i,p}\boldsymbol{X}_{j-i,T}$ for $j=p+1,...,T$.
\item[(2)] Define $\hat{\boldsymbol D}_T^*(v,\omega)$ in the same way as $\hat{\boldsymbol D}_T(v,\omega)$ but with the random variables $\boldsymbol{X}_{t,T}$ replaced by their bootstrap replicates $\boldsymbol{X}_{t,T}^*$.
\item[(3)] Set 
$\hat D_T^*:= \|[\sup\limits_{(v,\omega)\in[0,1]^2}|[\hat{\boldsymbol D}_T^*(v,\omega)]_{a,b}|]_{a,b=1,...,d} \|_F$
\end{itemize}
\end{algo}

To motivate this procedure note that for an increasing order $p$ the approximation in \eqref{arpfitt} becomes more accurate, which together with the consistency of the estimators $\hat a_{j,p}$ $(j=1,...,p)$ implies that the simulated data $\hat X_{t,T}^*$ and the original data $X_{t,T}$ have 'similar' distributional properties under the null hypothesis. We will show that this similarity carries over to the distributions of $\hat D_T^*$ and $\hat D_T$ under $H_0$, while, under $H_1$,  $\hat D_T^*$ still approximates the distribution of $\hat D_T$ under a stationary model with spectral density matrix $\int_0^1\boldsymbol{f}(u,\lambda)du$, whereas  $\hat D_T$ tends to some positive constant because of Theorem \ref{hauptsatz} b). In order to be able to state the necessary assumptions for this procedure to work  we introduce the process
\be
\boldsymbol{X}_{t}^{AR}(p):=\sum_{j=1}^p\boldsymbol{a}_{j,p}\boldsymbol{X}_{t-j}^{AR}(p)+\boldsymbol{Z}_t^{AR}(p), \label{arptrue}
\ee
where the $\boldsymbol{Z}_j^{AR}(p)$ are independent centered $\mathbb{R}^d$ valued Gaussian random vectors with covariance-matrix 
\bea
\boldsymbol{\Sigma}_p=\mathbb{E}\Big((\boldsymbol{X}_t-\sum_{j=1}^p\boldsymbol{a}_{j,p}\boldsymbol{X}_{t-j})(\boldsymbol{X}_t-\sum_{j=1}^p\boldsymbol{a}_{j,p}\boldsymbol{X}_{t-j})^T\Big)
\eea
[recall that $\boldsymbol{X}_t$ stands for a stationary process with spectral density matrix $\lambda \mapsto \int_0^1\boldsymbol{f}(u,\lambda) du$]. This means, that the process $\boldsymbol{X}_{t}^{AR}(p)$ corresponds to the 'best' approximation of the process \eqref{AR1} by an AR($p$) model. With this notation in mind we can state the  following Theorem which formally establishes the accuracy of the proposed bootstrap method.
\begin{satz}\label{theorembootstrap}
Suppose Assumption \ref{ass1} and \ref{annahmenbootstrap} are fulfilled.  Furthermore, let the following conditions on the growth rate of $p=p(T)$, the estimated AR parameters $\boldsymbol{\hat{a}}_{j,p}$ and its true counterparts  $\boldsymbol{{a}}_{j,p}$ and $\boldsymbol{{a}}_{j}$ be satisfied:
\begin{itemize}
\item[i)] There exist sequences $p_{min}(T)$ and $p_{max}(T)$ such that the order $p$ of the fitted autoregressive process satisfies $p=p(T) \in [p_{min}(T),p_{max}(T)]$ with $p_{max}(T)\geq p_{min}(T)\rightarrow \infty$ and
\be \label{pratebootstrap}
p^{3}_{max}(T)\sqrt{\log(T)/T}=O(1). 
\ee
\item[ii)] The estimators for the AR parameters satisfy
\be
\max_{1\leq j\leq p(T)}\|\boldsymbol{\hat{a}}_{j,p}-\boldsymbol{a}_{j,p}\|_{\infty}=O_P(\sqrt{\log(T)/T}) \label{ARparameterRate}
\ee
uniformly in $p$.
\item[iii)]The matrices $\boldsymbol{\Sigma}_p$ and $\boldsymbol{\Sigma}$ satisfy $\|\boldsymbol{\hat{\Sigma}}_{p}-\boldsymbol{\Sigma}\|_{\infty}\stackrel{P}{\longrightarrow}0$.
\end{itemize}
Then, as $T \rightarrow \infty$, we have 
$
\{\sqrt{T}\hat{\boldsymbol{ D}}_T^*(v,\omega)\}_{(v,\omega)\in[0,1]^2}\Rightarrow\{\boldsymbol{G}_{H_0}(v,\omega)\}_{(v,\omega)\in[0,1]^2}
$
conditionally on $\mathbb{X}(T)$,
where $\boldsymbol{G}_{H_0}(v,\omega)$ is a Gaussian process with covariance-kernel
\begin{align*}
&\text{Cov}([\boldsymbol{G}_{H_0}(v_1,\omega_1)]_{a_1,b_1},[\boldsymbol{G}_{H_0}(v_2,\omega_2)]_{a_2,b_2})\\
=&\frac{1}{2\pi} v_1v_2(\min(v_1,v_2)-v_1v_2)  \int_0^{\min(\omega_1,\omega_2)} \int_0^1f_{a_1b_2}(u,\lambda)du \int_0^1f_{b_1a_2}(u,-\lambda)du  d\lambda .
\end{align*}
\end{satz}

Note that the covariance kernel of the above limiting distribution is the 'same' as in Theorem \ref{hauptsatz} a) with the only difference that the spectral density $\boldsymbol{f}$ is replaced by the time-averaged version $\int_0^1 \boldsymbol{f}(u,\lambda) du$, which corresponds to   the [in the $L_2$-sense] best approximation of $\boldsymbol{f}(u,\lambda)$ by a stationary spectral density. Concerning the assumptions of Theorem \ref{theorembootstrap} we further note that they are rather standard in this framework and  are, for example, fulfilled for the Yule-Walker or the Least-Squres estimators [see \cite{bergpappolitis2010} or \cite{kreisspappol2011} for more details]. In order to obtain a decision rule which yields a level $\alpha$ test for the null hypothesis \eqref{null}, we now proceed according to the following algorithm:
\begin{algo}\label{algorithmtest}\mbox{}
\begin{itemize}
\item[1)] Calculate the test statistic (\ref{supstat}) using the observed data $\mathbb{X}(T)=\{\boldsymbol{X}_{1,T},...,\boldsymbol{X}_{T,T}\}$.
\item[2)] Choose $p \in \N$ and compute the estimates $(\boldsymbol{\hat{a}}_{1,p},...,\boldsymbol{\hat{a}}_{p,p},\boldsymbol{\hat\Sigma}_p)$ to fit an AR($p$) model to the observed data.
\item[3)] Apply Algorithm \ref{algo1} $B \in \N$ times and denote the resulting replicates of the test statistic by $\hat D_{T,i}^*$ for $i=1,...,B$.
\item[4)] Estimate the $(1-\alpha)$ quantile of (\ref{supstat}) by the corresponding empirical quantile of the sample $\{\hat D_{T,1}^*,...,\hat D_{T,T}^*\}$, i.e.  compute the order statistic $(\hat D_T^*)_{T,1},...,(\hat D_T^*)_{T,B}$ and define 
\be \label{bootstrapquantile} \hat q_{1-\alpha}:=(\hat D_T^*)_{T,\lfloor(1-\alpha)\rfloor B}.\ee
\item[5)] Reject $H_0$ if
\be \label{test} \|[\sup\limits_{(v,\omega)\in[0,1]^2}|[\hat{\boldsymbol D}_T(v,\omega)]_{a,b}|]_{a,b=1,...,d} \|_F >\hat q_{1-\alpha}.\ee
\end{itemize}
\end{algo}

Theorem \ref{hauptsatz} a) and Theorem \ref{theorembootstrap} yield that \eqref{test} corresponds to an asymptotic level $\alpha$ test for the null hypothesis of a constant dependency structure which is consistent because of Theorem \ref{hauptsatz} b) and Theorem \ref{theorembootstrap}. Note that by using Theorem \ref{hauptsatz}, \ref{theorembootstrap} and applying the continuous mapping theorem we can derive several other test statistics in order to validate the null hypothesis of stationarity, for instance  $\sup_{a,b} \sup\limits_{(v,\omega)\in[0,1]^2}|[\hat D_T(v,\omega)]_{a,b}|$ or  any  matrix norm of  $[\sup\limits_{(v,\omega)\in[0,1]^2}|[\hat D_T(v,\omega)]_{a,b}|]_{a,b=1,...,d}$. For the sake of brevity, we will restrict ourselves to the proposed Frobenius norm approach and demonstrate in Section \ref{simulationsection} that this yields very good results from a practical point of view.

\begin{rem} \label{gaussian} ~~ \\
{\rm
We noted earlier that the assumption of Gaussianity can be weakened, allowing our method to be applied to a much richer class of linear processes. We introduced this assumption merely to enhance the readability of the proofs of Theorems \ref{hauptsatz} and \ref{theorembootstrap}, and the respective statements can be extended to the more general setting of non normally distributed residuals by using the same main arguments but employing much more notation and a lot of additional calculations. Furthermore, for the bootstrap procedure to produce valid results in this case, one has to employ a different method for simulating the innovations $\boldsymbol{Z}_t^*$ in part 2) of Algorithm \ref{algo1}, which, in addition to mimicking  the second order properties of the underlying time series, must be able to imitate the fourth cumulant of the innovation process of $\mathbb{X}(T)$ [see \cite{krepap2012} for more details on this]. 
}
\end{rem}

\begin{rem} \label{identification} ~~ \\
{\rm 
If the testing procedure summarized in Algorithm \ref{algorithmtest} rejects the null hypothesis of a constant dependency structure for some $d$ variate time series  $\boldsymbol{X}_{t,T}=(X_{t,T,1},...,X_{t,T,d})$,   it might be of interest to further investigate the characteristics of the specific deviation from stationarity present in the data. More specifically it is sometimes  of interest  to determine whether, in fact, all $d$ components of the time series are non stationary or if there exists a subset $S:=\{i_1,...,i_{d'}\}\subset\{1,...,d\}$ of indices, such that the $d'$- variate time series   $\boldsymbol{X}_{t,T}'=(X_{t,T,i_1},...,X_{t,T,i_{d'}})$   can be still viewed as having a constant second-order structure. The construction of   such a refined analysis is enabled by part b) of Theorem \ref{hauptsatz}, and a natural way would be to proceed along the following steps: For each pair $(a,b)\in\{1,...,d\}^2$ we set 
\begin{equation} \label{identicomp2}
c(a,b):=
\begin{cases}
1\quad &\text{if}\quad\sqrt{T}\hat D_T(a,b)>\hat \epsilon_{T,a,b}\\
0\quad &\text{if}\quad \sqrt{T}\hat D_T(a,b)\leq \hat \epsilon_{T,a,b},
\end{cases}
\end{equation}
where $ \hat D_T(a,b):=\sup\limits_{v,\omega\in[0,1]}|[\boldsymbol{D}(v,\omega)]_{a,b}|$, $(a,b)\in\{1,...,d\}^2$ and $\hat \epsilon_{T,a,b}$ denotes some threshold sequence which tends to infinity but at a smaller rate than $\sqrt{T}$, as $T$ increases. We then compute 
$$d':=\max\{x\in\{1,...,d\}|\exists a_1,...,a_x\in\{1,...,d\}: c(a,b)=0\quad\forall (a,b)\in\{a_1,...,a_x\}\}$$
and choose elements  $i_1,...,i_{d'}$ with $c(i_a,i_b)=0$ for all $a,b=1,...,d'$ in order to construct the set $S$. Obviously this choice might not be unique, but in this case a repeated application of this procedure would result in not only one but maybe more subseries $\boldsymbol{X}_{t,T}'=(X_{t,T,i_1},...,X_{t,T,i_{d'}})$ which can be considered stationary.  Theorem \ref{hauptsatz} b) implies that, by using this algorithm, both the non-stationary and the stationary components will be detected as such with asymptotic probability one. Motivated by hard thresholding [see \cite{thresholding2}] we propose to use
\bea
\hat \epsilon_{T,a,b}=T^{\gamma} \sqrt{2 \hat V_{T,a,b} \log(d(d+1)/2)}
\eea
where $\gamma \in (0,1/2)$ serves as a tuning parameter and $\hat V_{T,a,b}=0.0125 \pi^{-2} T^{-2} \sum_{t=1}^T X_{t,T,a}^2 \sum_{t=1}^T X_{t,T,b}^2$ is a  proxy for the kind of variation in the process [the reasoning behind this proxy is obtained by maximizing the variance kernel in Theorem \ref{hauptsatz} over $v_i$ and $\omega_i$ in the White Noise case]. Note, however, that a detailed analysis of this problem is far too complex to be included in this paper as well.
}
\end{rem}

\section{Several Extensions of our approach}
\def\theequation{4.\arabic{equation}}
\setcounter{equation}{0}

The problem of testing for a constant dependency structure is, of course, only one issue which practitioners have an interest in when working with locally stationary processes. Although this question has drawn the largest attention in the literature  when it comes to goodness-of-fit testing in this framework [see the introduction for a brief overview over existing procedures], the problem of consistent parameter estimation in parametric or semiparametric models has received considerable interest as well [see \cite{dahlhaus1997}, \cite{dahlhaus2000} or \cite{dahlpolo2009} among many others]. Moreover, there exist several articles in which tests for a (semi-) parametric structure are developed [cf. \cite{Taniguchi}, \cite{sergpapa2009} or \cite{preusssemi}] or in which the spectral densities of different locally stationary time series models are compared with each other in order to obtain discriminating or clustering procedures [see \cite{sakitani2004}, \cite{huang2004}, \cite{chandler2006} or \cite{piotr2009}]. For the sake of brevity, we will not  revise all these procedures here [an overview can be found in \cite{dahlhaus2011}], but we intend to stress that in the framework of parametric estimation or testing essentially all existing approaches   either require the choice of at least a window length in the segmentation of the data or are based on the pre-periodogram which usually performs poorly in finite sample situations.  In the following we briefly indicate how our approach might be extended to these problems in order to develop methods for treating these various questions in a way, which does not depend on the choice of regularizing parameters.

\paragraph{Parameter Estimation:}

Suppose we want to fit a locally stationary model with  time-varying spectral density $\boldsymbol{f}_\theta(u,\lambda)$, $\theta \in \Theta \subset \R^p$ to the observations $\boldsymbol{X}_{1,T},...,\boldsymbol{X}_{T,T}$, where  $\Theta$ is assumed to be compact and $\boldsymbol{f} (u,\lambda)=\boldsymbol{f}_{\theta_0}(u,\lambda)$ holds for some parameter $\theta_0$ which lies in the interior of $\Theta$. By furthermore ensuring that $\theta_0$ is the unique minimzer of the Kolmogorov-Smirnov distance
\bea
KS(\boldsymbol{f},\boldsymbol{f}_\theta)=\min_{v,\omega}\frac{v}{2\pi} \Big \| \int_0^{\omega\pi}\int_0^v ( \boldsymbol{f}(u,\lambda)- \boldsymbol{f}_\theta(u,\lambda)) dud\lambda \Big \|_F
\eea

between $\boldsymbol{f}$ and the parametric spectral density $\boldsymbol{f}_\theta$, a natural estimator for $\theta_0$ is given by
\bea
\hat \theta_T=\text{argmin}_{\theta \in \Theta} \sup_{v,\omega}\Big \|  \frac{v}{T}\sum_{k=1}^{\lfloor \omega   \lfloor vT/2\rfloor   \rfloor}\boldsymbol I_{2\lfloor vT/2\rfloor}\big(\lambda_{k,2\lfloor vT/2\rfloor} \big )- \int_0^{\omega\pi}\int_0^v \boldsymbol{f}_\theta(u,\lambda)  \Big \|_F.
\eea

Here, the integrated spectral density is replaced by the corresponding sum over the  usual periodogram just as it is done in deriving the statistic $\hat D_T$ in Section \ref{Teststatsection}. Consistency of this estimator is easily derived by using the methods used in the proof of Theorem \ref{hauptsatz} while the derivation of the asymptotic distribution is more involved and therefore  subject to future research.

\paragraph{Testing for parametric hypothesis:} 

As mentioned in the above discussion, one question of interest is whether the process belongs to the considered parametric family at all, i.e. if there exists some $\theta_0$ such that  $\boldsymbol{f} (u,\lambda)=\boldsymbol{f}_{\theta_0}(u,\lambda)$ holds. One natural approach for this problem is to compare the integrated periodogram with the corresponding spectral distribution function using the estimator $\hat \theta_T$, by computing 
\bea
\sup_{v,\omega} \Big \|  \frac{v}{T}\sum_{k=1}^{\lfloor \omega   \lfloor vT/2\rfloor   \rfloor}\boldsymbol I_{2\lfloor vT/2\rfloor}\big(\lambda_{k,2\lfloor vT/2\rfloor} \big )- \int_0^{\omega\pi}\int_0^v \boldsymbol{f}_{\hat \theta_T}(u,\lambda)  \Big \|_F.
\eea
The estimation of the corresponding quantiles of this statistic can be conducted using bootstrap methods which create pseudo observations by using  the fitted locally stationary model. A detailed theoretical treatment, however, is beyond the scope of this paper.

\paragraph{Comparing spectral densities for discriminant or cluster analysis:} 

A third area of interest are methods for comparing two locally stationary time series models in order to test  for equality of the corresponding spectral densities or to  obtain discriminant and clustering algorithms. Such procedures can be once more obtained by estimating the Kolmogorov-Smirnov  difference between the corresponding spectral densities by using the usual periodogram as described in this paper. In order to derive a formal test for the equality of two time-varying spectral densities,  resampling methods are required again, while discriminant and cluster algorithms directly follow.

\section{Implementation and finite sample properties} \label{simulationsection}
\def\theequation{5.\arabic{equation}}
\setcounter{equation}{0}
In this section we present the results of a simulation study which investigates the finite sample properties of the proposed bootstrap based test \eqref{test} by analyzing how well the proposed procedure approximates the nominal level under   the null hypothesis and  how much the power is achieved for different alternatives. In the latter case, we also compare our new approach to existing procedures, and it will be demonstrated that it performs quite well even in a rather 'unfair' comparison in which we take the best obtained results [by varying the underlying regularization parameters] of competing procedures as reference values.  All reported empirical quantiles are based on 1000 simulation runs with each run containing 200 bootstrap replications.

For the implementation we choose the order of the AR-model, which is fitted for generating the bootstrap replicates according to Algorithm \ref{algo1}, as the minimizer of the AIC criterion dating back to \cite{akaike1973}, which becomes
\begin{align*} 
\hat p = \text{argmin}_p \frac{2\pi}{T} \sum_{k=1}^{{T}/{2}} \Big( \log(\text{det}[\boldsymbol{f}_{\hat \theta(p)}(\lambda_{k,T})])+\text{tr}[(\boldsymbol{f}_{\hat \theta(p)}(\lambda_{k,T}))^{-1}\boldsymbol I_T(\lambda_{k,T})] \Big) + p/T
\end{align*}
in the context of stationary processes due to \cite{whittle1}. Here, $\boldsymbol f_{\hat \theta(p)}$ denotes the spectral density of the fitted stationary AR($p$) process and $\boldsymbol I_T$ is the usual stationary periodogram with the corresponding Fourier frequencies $\lambda_{k,T}=2\pi k/T$.\\
Furthermore, in the simulation study and in the calculation of $p$-values in the data examples, we only considered data sets of size $T=2^i$ for some $i\in\N$ and values $v\in\{2^i T^{-1}|1\leq i\leq\lfloor\log_2(T/2)\rfloor\}$ for the computation of the test statistic $\hat D_T$ and its bootstrap replicates $\hat D_T^*$. This implies that the periodogram estimates $\boldsymbol{I}_{2\lfloor vT/2\rfloor}$ and $\boldsymbol{I}_{2\lfloor T/2\rfloor}$ are based on blocks of length equal to a power of 2, resulting in a significant reduction of computational runtime.
\paragraph{Size of the test:}
We first examine how well the nominal level is approximated under the null hypothesis. For this reason, we start by considering the univariate models
\begin{align}
X_{t,T}&=\theta Z_{t-1}+Z_t\label{null1}\\
X_{t,T}&=\phi X_{t-1,T}+Z_t\label{null2},
\end{align}
where  $Z_i$ denote independent standard normal distributed random variables, and continue with a size investigation for the bivariate models
\begin{align}
\boldsymbol{X}_{t,T}&=
\begin{pmatrix}
\theta & 0.2\\
0.2 &\theta
\end{pmatrix}
\boldsymbol{Z}_{t-1}+\boldsymbol{Z}_{t}\label{null3}\\
\boldsymbol{X}_{t,T}&=
\begin{pmatrix}
\phi& 0.2\\
0.2 &\phi
\end{pmatrix}
\boldsymbol{X}_{t-1,T}+\boldsymbol{Z}_{t}\label{null4},
\end{align}
where the $\boldsymbol{Z}_i$ denote independent bivariate normal distributed random vectors with unit covariance matrix. The estimated rejection probabilities of the test \eqref{test} for the models \eqref{null1} -- \eqref{null4} with different values of $\theta$ and $\phi$ are presented in Table \ref{table1} for the univariate case and Table \ref{table2} for the multivariate situation. In both cases our new test  seems to be conservative for $T=64$ while the approximation becomes better for increasing $T$, yielding satisfying results for a sample size of $256$ data.

\begin{table}[htbp]
\centering
\begin{tabular}{|c|c|c|c|c|c|c|c|c|c|c|}
\hline
\multicolumn{ 1}{|c|}{}&\multicolumn{ 6}{|c|}{$\text{H}_0$: Model \eqref{null1}}&\multicolumn{ 4}{|c|}{ $\text{H}_0$: Model \eqref{null2} } \\ \hline
 & \multicolumn{ 2}{c|}{$\theta=0.5$}& \multicolumn{ 2}{c|}{$\theta=0$} & \multicolumn{ 2}{c|}{$\theta=-0.5$} & \multicolumn{ 2}{c|}{$\phi=0.5$} & \multicolumn{ 2}{c|}{$\phi=-0.5$} \\ \hline
$T$ & 5 \% & 10 \% &5 \% & 10 \% & 5 \% & 10 \% & 5 \% & 10 \% & 5 \% & 10 \% \\ \hline
64  & 0.031 & 0.070 & 0.036 & 0.078 & 0.043 & 0.095 & 0.044 & 0.094 & 0.023 & 0.068 \\ \hline
128  & 0.044 & 0.089 & 0.039 & 0.097 & 0.033 & 0.069 & 0.044 & 0.101 & 0.044 & 0.097 \\ \hline
256  & 0.043 & 0.094 & 0.052 & 0.098 & 0.038 & 0.077 & 0.038 & 0.082 & 0.045 & 0.105 \\ \hline
\end{tabular}
\caption{{\it Empirical rejection frequencies of the test \eqref{test} for the models \eqref{null1} and \eqref{null2} with different choices for the parameters $\theta$, $\phi$ and the sample size $T$ at nominal levels $\alpha=0.05$ and $\alpha=0.1$.}}
\label{table1}
\end{table} 

\begin{table}[htbp]
\centering
\begin{tabular}{|c|c|c|c|c|c|c|c|c|}
\hline
\multicolumn{ 1}{|c|}{}&\multicolumn{ 4}{|c|}{$\text{H}_0$: Model \eqref{null3}}&\multicolumn{ 4}{|c|}{ $\text{H}_0$: Model \eqref{null4} } \\ \hline
 & \multicolumn{ 2}{c|}{$\theta=0.5$} & \multicolumn{ 2}{c|}{$\theta=-0.5$} & \multicolumn{ 2}{c|}{$\phi=0.5$} & \multicolumn{ 2}{c|}{$\phi=-0.5$} \\ \hline
$T$ & 5 \% & 10 \% &5 \% & 10 \% & 5 \% & 10 \% & 5 \% & 10 \% \\ \hline
64  &  0.017 & 0.039 & 0.010 & 0.041 & 0.034 & 0.094 & 0.011 & 0.034 \\ \hline
128 & 0.030 & 0.061 & 0.018 & 0.051 & 0.038 & 0.091 & 0.020 & 0.057 \\ \hline
256 & 0.037 & 0.080 & 0.039 & 0.079 & 0.039 & 0.093 & 0.038 & 0.085 \\ \hline
\end{tabular}
\caption{{\it Empirical rejection frequencies of the test \eqref{test} for the models \eqref{null3} and \eqref{null4} with different choices for the parameters $\theta$, $\phi$ and the sample size $T$ at nominal levels $\alpha=0.05$ and $\alpha=0.1$.}}
\label{table2}
\end{table} 

\paragraph{Power of the test:}
We now continue by   investigating the power of the test \eqref{test} in various non stationary time series models. For this reason, we simulated data sets $\mathbb{X}(T)$ according to the following set of univariate and multivariate locally stationary processes
\begin{align}
X_{t,T}&=(1+\frac{t}{T})Z_t\label{alt1}\\
X_{t,T}&=-0.9\sqrt{t/T}X_{t-1,T}+Z_t\label{alt2}\\
X_{t,T}&=
\begin{cases}
0.5X_{t-1,T}+Z_t &\text{if}\quad1\leq t\leq\frac{T}{2}\\
-0.5X_{t-1,T}+Z_t &\text{if}\quad\frac{T}{2}+1\leq t\leq T
\end{cases}\label{alt3}\\
\boldsymbol{X}_{t,T}&=
1.4t/T \boldsymbol{A} \boldsymbol{X}_{t-1,T}+\boldsymbol{Z}_t, \label{alt5}\\
\boldsymbol{X}_{t,T}&=\sin(2\pi t/T)
\boldsymbol{A}\boldsymbol{X}_{t-1,T}+\boldsymbol{Z}_t, \label{alt6} \\
\boldsymbol{X}_{t,T}&=
\boldsymbol{A}\boldsymbol{X}_{t-1,T}+2\sin(2\pi t/T)\boldsymbol{Z}_t, \label{alt7}
\end{align}
with
\bea
\boldsymbol{A}=\begin{pmatrix}
0.6 & 0.2\\
0 & 0.3
\end{pmatrix}
 \text{, } \quad
\boldsymbol{\Sigma}=\begin{pmatrix}
1 & 0.3\\
0.3 & 1 
\end{pmatrix},
\eea
 and where in each model $\{Z_t\}_{t=1,...,T}$ and  $\{\boldsymbol{Z}_t\}_{t=1,...,T}$ denote  sequences of univariate/bivariate mean zero   Gaussian distributed random variables   with the variance of $Z_t$ being equal to one while we assume that $\boldsymbol{Z}_t$ has the (co)variance matrix $\boldsymbol{\Sigma}$. The processes \eqref{alt6} and \eqref{alt7} were considered in the simulation study of \cite{rao2012} which we would like to compare our testing procedure with. Note that the process \eqref{alt3} is strictly speaking not inluded in our theoretical framework, but since we conjecture that our approach can be extended to include break-point models, the investigation of the 
 corresponding finite sample properties is of interest as well. 

The obtained estimates for the rejection probabilities of the test \eqref{test} are sumarized in Table \ref{table3}  for the models \eqref{alt1} -- \eqref{alt6}, from which we observe that, for all considered alternatives, the rejection frequencies are far above the nominal levels and increasing with the sample size, just as the theoretical results indicate. In order to compare the performance  of our new procedure with existing approaches in the univariate case, we make use of the comprehensive comparison study   in \cite{detprevet2011b}. In this paper, two different estimators of the alternative Kolmogorov-Smirnov distance \eqref{Distance} were used to derive two tests for stationarity, whose finite sample properties  were then compared to those obtained from the approaches of \cite{paparoditis2010}, \cite{rao2010} and \cite{detprevet2010} by using the models \eqref{alt1}--\eqref{alt3}. The two different test statistics proposed in \cite{detprevet2011b} estimate  \eqref{Distance} either by a summarized local periodogram [meaning that a window-length has to be chosen] or by a Riemann-sum over the pre-periodogram of \cite{neumsach1997}. In a comprehensive simulation study it was found out that the local periodogram approach works much better than employing the pre-periodogram  and that it also performs quite well compared to the alternative tests included in the investigation. In Table \ref{table4} we  partly restate the results of the simulation study in \cite{detprevet2011b} by presenting  intervals  containing both the lowest and the highest rejection frequency obtained for the models \eqref{alt1}--\eqref{alt3} at  nominal level $\alpha=0.1$ [note that the wide range of the intervals is due to several tests included in the investigation which additionally all depend on several regularization parameters yielding different results]. In addition we highlight  the results for   models, in which the local periodogram Kolmogorov-Smirnov approach performed best for some window-length by printing them boldly. We see that this was the case for the processes \eqref{alt1} and \eqref{alt2}, and that, for these two models, our new approach performs quite well compared to existing ones. Note that, for   model \eqref{alt2}, the old local periodogram procedure works slightly better, but only if the window-length is chosen 'correctly'. In fact, the rejection frequencies can heavily depend on the choice of this regularization parameter, and  it is therefore advisable to employ our new procedure in these cases. 

Concerning the process \eqref{alt3}, it turns out that the approach of \cite{detprevet2010} has the largest power, but again only for one specific choice of the window-length. In fact,   by choosing a different value for this regularization parameter, the stated best rejection frequency of $0.922$ for $T=256$ can go down to $0.654$, implying that our new procedure would yield higher rejection frequencies even in this case. So, we conclude  as a summary for the one dimensional case that, by using the new approach, the practitioner does not face the possibility of obtaining artificial results by varying the underlying regularization parameters, resulting in a much more reliable procedure to test the null hypothesis of a constant dependency structure, which additionally performs reasonably well compared to existing approaches [note that we compared our method with a large number of competing approaches for several choices of the corresponding tuning parameters].

For the multivariate models \eqref{alt5} -- \eqref{alt7} we compare the performance of our procedure with the proposed test of \cite{rao2012} in which we choose all four tuning parameters as proposed by these authors and use their non-bootstrap procedure since this  yields higher power estimates as discussed in their simulation study. The results of this approach are summarized in Table \ref{table5} and if we compare the rejection frequencies with the ones of our test displayed in Table \ref{table3}, we see that our procedure clearly outperforms the test of \cite{rao2012} for the processes \eqref{alt5} and \eqref{alt6} while the opposite is the case when data come from the model \eqref{alt7}. It is, however, a well known fact that Kolmogorov-Smirnov approaches lack power in periodic time series models, so it is  already remarkable that our method performs much better for the process \eqref{alt6} [in fact, for a large number of considered non-periodic  parameter functions, our procedure has a much larger power]. Thus, to summarize the multivariate setting, our decision rule \eqref{test} is (at least) quite competitive to the existing approach of \cite{rao2012} and requires zero instead of three tuning parameters in the calculation of the test statistic.

\begin{table}[htbp]
\centering
\begin{tabular}{|c|c|c|c|c|c|c| c|c|c|c|c|c|}
\hline
\multicolumn{ 1}{|c|}{}&\multicolumn{ 2}{|c|}{ Model \eqref{alt1}}&\multicolumn{ 2}{|c|}{Model \eqref{alt2}}&\multicolumn{ 2}{|c|}{Model \eqref{alt3}}  & \multicolumn{ 2}{|c|}{Model \eqref{alt5}} & \multicolumn{ 2}{|c|}{Model \eqref{alt6}} & \multicolumn{ 2}{|c|}{Model \eqref{alt7}} \\ \hline
$T$   &5 \% & 10 \% & 5 \% & 10 \% & 5 \% & 10 \%& 5 \% & 10 \%& 5 \% & 10 \% & 5 \% & 10 \% \\ \hline
64  & 0.341 & 0.474 & 0.088 & 0.256 & 0.096 & 0.214  & 0.100 & 0.218 & 0.070 & 0.194 & 0.060 & 0.150 \\ \hline
128  & 0.698 & 0.788 & 0.311 & 0.501 & 0.182 & 0.368   & 0.292 & 0.462 & 0.150 & 0.308 & 0.084 & 0.174 \\ \hline
256  & 0.958 & 0.971 & 0.737 & 0.829 & 0.416 & 0.727   & 0.684 & 0.796  &0.370 & 0.608 & 0.096 & 0.186 \\ \hline
\end{tabular}
\caption{{\it  Empirical rejection frequencies of the test \eqref{test} for the models \eqref{alt1} -- \eqref{alt6}, sample sizes $T\in\{64,128,256\}$ and nominal levels $\alpha=0.05$, $\alpha=0.1$.}}
\label{table3}
\end{table}

\begin{table}[htbp]
\centering
\begin{tabular}{|c|c|c|c|}
\hline
\multicolumn{ 1}{|c|}{T}&\multicolumn{ 1}{|c|}{ Model \eqref{alt1} at 10 \%}&\multicolumn{ 1}{|c|}{Model \eqref{alt2} at 10 \%}&\multicolumn{ 1}{|c|}{Model \eqref{alt3} at 10 \%}   \\ \hline
64  & [0.126,{\bf 0.444}] & [0.100,{\bf 0.328}] & [0.056,0.344]  \\ \hline
128  & [0.16,{\bf 0.772}] & [0.114,{\bf 0.578}] & [0.116,0.566]   \\ \hline
256  & [0.226,{\bf 0.978}] & [0.210,{\bf 0.868}] & [0.176,0.922]  \\ \hline 
\end{tabular}
\caption{{\it  Empirical rejection frequencies for the tests  of \cite{paparoditis2010}, \cite{rao2010}, \cite{detprevet2010} and \cite{detprevet2011b} for the models \eqref{alt1} -- \eqref{alt3}, sample sizes $T\in\{64,128,256\}$ and nominal level   $\alpha=0.1$. The bold results correspond to the cases, where the local periodogram approach of \cite{detprevet2011b} works best for some specific window-length.}}
\label{table4}
\end{table}

\begin{table}[htbp]
\centering
\begin{tabular}{|c|c|c|c|c|c|c|}
\hline
\multicolumn{ 1}{|c|}{T}&\multicolumn{ 2}{|c|}{ Model \eqref{alt5}}&\multicolumn{ 2}{|c|}{Model \eqref{alt6}}&\multicolumn{ 2}{|c|}{Model \eqref{alt7}}   \\ \hline
$T$   &5 \% & 10 \% & 5 \% & 10 \% & 5 \% & 10 \%  \\ \hline
64  & 0.051 & 0.098 & 0.050 & 0.132 & 0.201 & 0.256   \\ \hline
128  & 0.055 & 0.122 & 0.104 & 0.214 & 0.551 & 0.699    \\ \hline
256   & 0.118 & 0.147 & 0.188 & 0.292 & 0.998 & 1.00    \\ \hline 
\end{tabular}
\caption{{\it  Empirical rejection frequencies for the test  of \cite{rao2012} for the models \eqref{alt5} -- \eqref{alt7}, sample sizes $T\in\{64,128,256\}$ and nominal levels $\alpha=0.05$, $\alpha=0.1$.}}
\label{table5}
\end{table} 

 {\bf Acknowledgements.}
This work has been supported in part by the
Collaborative Research Center ``Statistical modeling of nonlinear
dynamic processes'' (SFB 823, Teilprojekt C1) of the German Research Foundation (DFG).

 %BEWEISE
\section{Appendix}
\def\theequation{6.\arabic{equation}}
\setcounter{equation}{0}

\subsection{Proof of Theorem \ref{hauptsatz}}
For notational convenience and without loss of generality we restrict ourselves to the case $d=1$, since the more general case is treated completely analogously using linearity arguments and the independence of the components of $\boldsymbol{Z}_t$. Throughout this chapter $C$ denotes some universal positive constant that may vary from line to line and $\epsilon'>0$ denotes some constant which can be arbitrarily small. \\

{\bf Proof of part a):} Following Theorems 1.5.4 and 1.5.7. in \cite{wellnervandervaart}, it suffices to show 
\begin{itemize}
\item[(1)] \textit{[Convergence of the finite dimensional distributions]} For every  $k\geq 1$, $(v_1,\omega_1),...,(v_k,\omega_k)\in[0,1]^2$ it holds
\begin{equation}\label{fidis}
\sqrt{T}\big( \hat D_{T}(v_1,\omega_1),...,\hat D_{T}(v_k,\omega_k)\big)\Rightarrow \big(G(v_1,\omega_1),...,G(v_k,\omega_k)\big),
\end{equation}
\item[(2)] \textit{[Stochastic Equicontinuity]}  The process $\{\sqrt{T}\hat{ D}_{T}(v,\omega)\}_{(v,\omega)\in[0,1]^2}$ is asymptotically stochastically equicontinuous, i.e. for every $\eta, \epsilon>0$ there exists a $\delta>0$ such that
\begin{align}
\label{Gleichstetigkeit}
\lim_{T\rightarrow \infty} P\Big( \sup_{y_1,y_2 \in [0,1]^2: d(y_1,y_2) <\delta } \sqrt{T}|\hat{ D}_{T}(y_1)-\hat{  D}_{T}(y_2)|> \eta \Big)< \epsilon ,
\end{align} 
where $d(y_1,y_2)$ denotes the distance measure $(|v_1-v_2|+|\omega_1-\omega_2|)^{\beta/2}$ between $y_1:=(v_1,\omega_1)$ and $y_2:=(v_2,\omega_2)$ for some fixed $0<\beta<1/3$.
\end{itemize}

\textbf{Proof of \eqref{fidis}:} Convergence of the finite dimensional distributions of the empirical process \linebreak $\{\sqrt{T}\hat{ D}_{T}(v,\omega)\}_{(v,\omega)\in[0,1]^2}$ can be obtained by showing that all cumulants of the random vector \linebreak{$[\sqrt{T}\hat D_{T}(v_i,\omega_i)]_{i=1,...,k}$} converge to the respective cumulants of the vector $[G(v_i,\omega_i)]_{i=1,...,k}$. Therefore we proceed by showing that the following claims hold under the null hypothesis:
\begin{itemize}
\item[(i)]$\E(\sqrt{T}\hat D_{T}(v,\omega))=o(1)$ for all $(v,\omega)\in[0,1]^2$.
\item[(ii)]$\text{Cov}(\sqrt{T}\hat D_{T}(v_1,\omega_1),\sqrt{T}\hat D_{T}(v_2,\omega_2))=\text{Cov}(G(v_1,\omega_1),G(v_2,\omega_2))+o(1)$ for all $(v_1,\omega_1),(v_2,\omega_2)\in[0,1]^2$.
\item[(iii)]$\text{cum}(\sqrt{T}\hat D_{T}(v_1,\omega_1),...,\sqrt{T}\hat D_{T}(v_l,\omega_l))=o(1)$ for all $l\geq 3$, $(v_1,\omega_1),...,(v_l,\omega_l)\in[0,1]^2$.
\end{itemize}
\textit{Proof of part (i):} We prove that the slightly more general statement
\begin{align}\label{expectation1}
\E\Big(\sqrt{T}(\hat D_T(v,\omega)-D(v,\omega))\Big)=o(1)
\end{align}
holds for all $(v,\omega)\in[0,1]^2$ both under $H_0$ and $H_1$, since it will reduce the arguments required in the proof  of part b)  significantly later on. It is straightfoward to see that \eqref{expectation1} is a consequence of 
\begin{align}\label{expectation2}
\E\Big(\sqrt{T}(\frac{v}{T}\sum_{k=1}^{\lfloor\omega\lfloor\frac{vT}{2}\rfloor\rfloor}I_{2\lfloor\frac{vT}{2}\rfloor}(\lambda_{k,2\lfloor\frac{vT}{2}\rfloor})-\frac{v}{2\pi}\int_0^{\omega\pi}\int_0^{v}f(u,\lambda)dud\lambda)\Big)=o(1),
\end{align}
and we thus focus on a proof of \eqref{expectation2} in the following. By using the notation $\psi_{t,T,l}=\boldsymbol{\Psi}_{t,T,l}$, $\psi_l(t/T):=\boldsymbol{\Psi}_l(t/T)$ and employing the locally stationary representation \eqref{defprocess1} of the time series $\mathbb{X}(T)$, we obtain
\begin{align}
&\E\big(\frac{v}{T}\sum_{k=1}^{\lfloor\omega\lfloor\frac{vT}{2}\rfloor\rfloor}I_{2\lfloor\frac{vT}{2}\rfloor}(\lambda_{k,2\lfloor\frac{vT}{2}\rfloor})\big)=\E\big(\frac{v}{T}\sum_{k=1}^{\lfloor\omega\lfloor\frac{vT}{2}\rfloor\rfloor}\frac{1}{4\pi\lfloor\frac{vT}{2}\rfloor}\sum_{p,q=0}^{2\lfloor\frac{vT}{2}\rfloor-1}X_{1+p}X_{1+q}\exp(-i\lambda_{k,2\lfloor\frac{vT}{2}\rfloor}(p-q))\big)\nonumber\\
=&\frac{v}{T}\sum_{k=1}^{\lfloor\omega\lfloor\frac{vT}{2}\rfloor\rfloor}\frac{1}{4\pi\lfloor\frac{vT}{2}\rfloor}\sum_{l,m=0}^{\infty}\sum_{p,q=0}^{2\lfloor\frac{vT}{2}\rfloor-1}\psi_{1+p,T,l}\psi_{1+q,T,m}\E\big(Z_{1+p-l}Z_{1+q-m}\big)\exp(-i\lambda_{k,2\lfloor\frac{vT}{2}\rfloor}(p-q)).\nonumber %\label{expectation3},
\end{align}
By the approximating property \eqref{apprbed} we get that this term equals
\begin{align}
\frac{v}{T}\sum_{k=1}^{\lfloor\omega\lfloor\frac{vT}{2}\rfloor\rfloor}\frac{1}{4\pi\lfloor\frac{vT}{2}\rfloor}\sum_{l,m=0}^{\infty}\sum_{p,q=0}^{2\lfloor\frac{vT}{2}\rfloor-1}\psi_l(\frac{1+p}{T})\psi_m(\frac{1+q}{T})\E\big(Z_{1+p-l}Z_{1+q-m}\big)\exp(-i\lambda_{k,2\lfloor\frac{vT}{2}\rfloor}(p-q))+O(\frac{1}{T})\label{expectation3},
\end{align}
i.e. the time varying coefficents $\psi_{t,T,l}$ can be replaced by the approximating functions $\psi_l(t/T)$ by making an error of order $O(T^{-1})$. The condition $\E(Z_iZ_j)=\delta_{ij}$ implies that the restriction $q=p-l+m$ has to hold, which yields that \eqref{expectation3} is equal to
\begin{align}
\frac{v}{T}\sum_{k=1}^{\lfloor\omega\lfloor\frac{vT}{2}\rfloor\rfloor}\frac{1}{2\pi}\sum_{l,m=0}^{\infty}\exp(-i\lambda_{k,2\lfloor\frac{vT}{2}\rfloor}(l-m))\frac{1}{2\lfloor\frac{vT}{2}\rfloor}\sum_{\substack{p=0\\0\leq p-l+m\leq 2\lfloor\frac{vT}{2}\rfloor-1}}^{2\lfloor\frac{vT}{2}\rfloor-1}\psi_l(\frac{1+p}{T})\psi_m(\frac{1+p-l+m}{T})+O(\frac{1}{T})\label{expectation4}.
\end{align}
Because of
\bea
&& \sum_{l,m=0}^{\infty}\Big | \sum_{\substack{p=0\\0\leq p-l+m\leq 2\lfloor\frac{vT}{2}\rfloor-1}}^{2\lfloor\frac{vT}{2}\rfloor-1}\psi_l(\frac{1+p}{T})\psi_m(\frac{1+p-l+m}{T})-\sum_{\substack{p=0}}^{2\lfloor\frac{vT}{2}\rfloor-1}\psi_l(\frac{1+p}{T})\psi_m(\frac{1+p-l+m}{T}) \Big | \\ && \leq  \sum_{l,m=0}^{\infty} |l-m| \sup\limits_{u\in[0,1]}|\psi_l(u)| \sup\limits_{u\in[0,1]}|\psi_m(u)|
\eea
the property \eqref{summ1} now shows that the restriction in the summation over $p$ in \eqref{expectation4} can be dropped by making an error of order $O(T^{-1})$. We thus get that \eqref{expectation4} is equal to
\begin{align*}
&\frac{v}{T}\sum_{k=1}^{\lfloor\omega\lfloor\frac{vT}{2}\rfloor\rfloor}\frac{1}{2\pi}\sum_{l,m=0}^{\infty}\exp(-i\lambda_{k,2\lfloor\frac{vT}{2}\rfloor}(l-m))\frac{1}{2\lfloor\frac{vT}{2}\rfloor}\sum_{p=0}^{2\lfloor\frac{vT}{2}\rfloor-1}\big[\psi_l(\frac{1+p}{T})\psi_m(\frac{1+p}{T})\\
&+\psi_l(\frac{1+p}{T})\big(\psi_m(\frac{1+p-l+m}{T})-\psi_m(\frac{1+p}{T})\big)\big]+O(\frac{1}{T})=E_{T,1}+E_{T,2}+O(\frac{1}{T}) ,
\end{align*}
with $E_{T,1}$ and $E_{T,2}$ being defined implicitly. The claim \eqref{expectation2}  follows if we show:
\begin{itemize}
\item[1)]$E_{T,1}=\frac{v}{2\pi}\int_0^{\omega\pi}\int_0^{v}f(u,\lambda)dud\lambda+O(\frac{1}{T})$.
\item[2)]$E_{T,1}=O(\frac{1}{T})$.
\end{itemize}
Regarding 1) we obtain that $E_{T,1}$ is equal to
\begin{align*}
&\frac{v}{T}\sum_{k=1}^{\lfloor\omega\lfloor\frac{vT}{2}\rfloor\rfloor}\frac{1}{2\pi}\sum_{l,m=0}^{\infty}\exp(-i\lambda_{k,2\lfloor\frac{vT}{2}\rfloor}(l-m))\frac{1}{v}\int_0^v\psi_l(u)\psi_m(u)du+O(\frac{1}{T})\\
=&\frac{\lfloor\frac{vT}{2}\rfloor}{T}\frac{1}{\lfloor\frac{vT}{2}\rfloor}\sum_{k=1}^{\lfloor\omega\lfloor\frac{vT}{2}\rfloor\rfloor}\int_0^vf(u,\lambda_{k,2\lfloor\frac{vT}{2}\rfloor})du+O(\frac{1}{T})=\frac{v}{2\pi}\int_0^{\omega\pi}\int_0^{v}f(u,\lambda)dud\lambda+O(\frac{1}{T}).
\end{align*}
With respect to 2) a Taylor expansion combined with the properties \eqref{summ1} and \eqref{summ2} implies that $|E_{T,2}|$ is bounded by
\begin{align*}
C\sum_{l,m=0}^{\infty}\sup\limits_{u\in[0,1]}|\psi_l(u)|\frac{|m-l|}{T}\sup\limits_{u\in[0,1]}|\psi_m'(u)|=O(\frac{1}{T}),
\end{align*}
which completes the proof of \eqref{expectation2}.\\

\textit{Proof of (ii):} For the calculation of the covariances we define the process
\be
\hat D_{T,1}(v,\omega):=\frac{v}{T}\sum_{k=1}^{\lfloor \omega\lfloor vT/2\rfloor\rfloor}I_{2\lfloor vT/2\rfloor}(\lambda_{k,2\lfloor vT/2\rfloor})
\ee
for $(v,\omega)\in[0,1]^2$ and note that it is $\hat D_T(v,\omega)=\hat D_{T,1}(v,\omega)-v^2 \hat D_{T,1}(1,\omega)$ [i.e. it is sufficient to calculate the covariance kernel of $\hat D_{T,1}(v,\omega)$]. Without loss of generality we restrict ourselves to the case $v_1\leq v_2$ and obtain 
\begin{align}
&\text{Cov}\big(\sqrt{T}\hat D_{T,1}(v_1,\omega_1),\sqrt{T}\hat D_{T,1}(v_2,\omega_2)\big)\nonumber\\
=&\frac{v_1v_2}{T}\sum_{k_1=1}^{\lfloor\omega_1\lfloor v_1T/2\rfloor\rfloor}\sum_{k_2=1}^{\lfloor\omega_2\lfloor v_2T/2\rfloor\rfloor}\text{cum}\big(I_{2\lfloor \frac{v_1T}{2}\rfloor}(\lambda_{k_1,2\frac{\lfloor v_1T}{2}\rfloor}),I_{2\lfloor \frac{v_2T}{2}\rfloor}(\lambda_{k_2,2\lfloor \frac{v_2T}{2}\rfloor})\big).\nonumber
\end{align}
By inserting the corresponding definitions we obtain that the above expression equals
\begin{align}
&\frac{v_1v_2}{T}\sum_{k_1=1}^{\lfloor\omega_1\lfloor v_1T/2\rfloor\rfloor}\sum_{k_2=1}^{\lfloor\omega_2\lfloor v_2T/2\rfloor\rfloor}\frac{1}{(2\pi)^24\lfloor\frac{v_1T}{2}\rfloor\lfloor\frac{v_2T}{2}\rfloor}\sum_{p_1,q_1=0}^{2\lfloor v_1T/2\rfloor-1}\sum_{p_2,q_2=0}^{2\lfloor v_2T/2\rfloor-1}\exp(-i\lambda_{k_1,2\lfloor\frac{v_1T}{2}\rfloor}(p_1-q_1))\nonumber\\
\times&\exp(-i\lambda_{k_2,2\lfloor\frac{v_2T}{2}\rfloor}(p_2-q_2))\sum_{l,m,n,o=0}^{\infty}\psi_l\psi_m\psi_n\psi_o\text{cum}(Z_{1+p_1-l}Z_{1+q_1-m},Z_{1+p_2-n}Z_{1+q_2-o})\label{vorProd},
\end{align}
where $\psi_l:=\psi_l(u)$ [note that we are under $H_0$].  
By means of the product theorem for cumulants [cf. Theorem 2.3.2 in \cite{brillinger1981}] we obtain
\begin{align*}
\text{cum}(Z_{a}Z_{b},Z_{c}Z_{d})=\text{cum}(Z_{a},Z_{d})\text{cum}(Z_{b},Z_{c})+\text{cum}(Z_{a},Z_{c})\text{cum}(Z_{b},Z_{d}),
\end{align*}
which implies that the calculation of the covariances can be split up into two parts $V_{T,1}$, $V_{T,2}$. We start by considering the first one and afterwars show that $V_{T,2}$ vanishes in the limit [i.e. $V_{T,1}$ is the dominating term]. Concerning $V_{T,1}$, we obtain with  $\text{cum}(Z_i,Z_j)=\delta_{ij}$  that the relations $q_2=p_1-l+o$ and $q_1=p_2-n+m$ have to hold, implying that $V_{T,1}$ equals
\begin{align}
&\frac{v_1v_2}{T}\sum_{k_1=1}^{\lfloor\omega_1\lfloor v_1T/2\rfloor\rfloor}\sum_{k_2=1}^{\lfloor\omega_2\lfloor v_2T/2\rfloor\rfloor}\frac{1}{(2\pi)^2}\sum_{l,m,n,o=0}^{\infty}\psi_l\psi_m\psi_n\psi_o\exp(-i\lambda_{k_1,2\lfloor\frac{v_1T}{2}\rfloor}(m-n))\exp(-i\lambda_{k_2,2\lfloor\frac{v_2T}{2}\rfloor}(l-o))\nonumber\\
\times&\frac{1}{4\lfloor\frac{v_1T}{2}\rfloor\lfloor\frac{v_2T}{2}\rfloor}\sum_{\substack{p_1=0\\0\leq p_1-l+o\leq2\lfloor\frac{v_2T}{2}\rfloor-1}}^{2\lfloor v_1T/2\rfloor-1}\sum_{\substack{p_2=0\\0\leq p_2-n+m\leq2\lfloor\frac{v_1T}{2}\rfloor}}^{2\lfloor v_2T/2\rfloor-1}\exp(-i(\lambda_{k_1,2\lfloor\frac{v_1T}{2}\rfloor}-\lambda_{k_2,2\lfloor\frac{v_2T}{2}\rfloor})(p_1-p_2)).\label{cov1}
\end{align}
We now show that we can drop the $l$, $m$, $n$ and $o$ terms in the restriction with the summation over $p_1$ and $p_2$ by making an error that is of order $O(1/T^{1-\epsilon'})$ for any $\epsilon'>0$. It is easy to see that the error which is made by changing the restriction in the summation over $p_1$ to $0\leq p_1\leq2\lfloor v_2T/2\rfloor-1$ is bounded by
\begin{align*}
&\frac{Cv_1v_2}{T}\sum_{k_1=1}^{\lfloor\omega_1\lfloor v_1T/2\rfloor\rfloor}\sum_{k_2=1}^{\lfloor\omega_2\lfloor v_2T/2\rfloor\rfloor}\sum_{l,m,n,o=0}^{\infty}|o-l|\psi_l\psi_m\psi_n\psi_o
\frac{1}{4\lfloor v_1T/2\rfloor\lfloor v_2T/2\rfloor}\Big|\sum_{\substack{p_2=0\\0\leq p_2-n+m\leq2\lfloor\frac{v_1T}{2}\rfloor}}^{2\lfloor v_2T/2\rfloor-1}\nonumber\\
\times&\exp(-i(\lambda_{k_1,2\lfloor v_1T/2\rfloor}-\lambda_{k_2,2\lfloor v_2T/2\rfloor})(-p_2))\Big|=:E_{T,\eqref{res1}}+E_{T,\eqref{res2}},
\end{align*}
where $E_{T,\eqref{res1}}$ and $E_{T,\eqref{res2}}$ denote the summations with respect to the conditions
\begin{align}
&|2k_1\lfloor v_2T/2\rfloor-2k_2\lfloor v_1T/2\rfloor|\leq2\lfloor v_1T/2\rfloor,\label{res1}\\
2\lfloor v_1T/2\rfloor<&|2k_1\lfloor v_2T/2\rfloor-2k_2\lfloor v_1T/2\rfloor|\leq2\lfloor v_1T/2\rfloor\lfloor v_2T/2\rfloor.\label{res2}
\end{align}
As for fixed $k_1$ there exists at most three value for $k_2$ such that \eqref{res1} holds it follows from \eqref{summ2}, that $E_{T,\eqref{res1}}=O(1/T)$. For the treatment of $E_{T,\eqref{res2}}$ we note that the geometric series formula and the identity $|1-\exp(ix)|=2|\sin(x/2)|$ yield
\begin{align*}
&\Big|\sum_{p_2=0}^{2\lfloor\frac{v_2T}{2}\rfloor-1}\exp(-i(\lambda_{k_1,2\lfloor\frac{v_1T}{2}\rfloor}-\lambda_{k_2,2\lfloor\frac{v_2T}{2}\rfloor})p_2)\Big|=\frac{|\big(1-\exp(-i(\lambda_{k_1,2\lfloor\frac{v_1T}{2}\rfloor}-\lambda_{k_2,2\lfloor\frac{v_2T}{2}\rfloor})2\lfloor\frac{v_2T}{2}\rfloor)\big)|}{|1-\exp(-i(\lambda_{k_1,2\lfloor\frac{v_1T}{2}\rfloor}-\lambda_{k_2,2\lfloor\frac{v_2T}{2}\rfloor}))|} \nonumber\\
\leq& \sin^{-1}\Big((\lambda_{k_1,2\lfloor\frac{v_1T}{2}\rfloor}-\lambda_{k_2,2\lfloor\frac{v_2T}{2}\rfloor})/2\Big).
\end{align*}
This property together with the bound $\sin(\pi x)\geq Cx$ for $x\in[0,1/2]$ implies that $E_{T,\eqref{res2}}$ is bounded by
\begin{align*}
&\frac{Cv_1v_2}{T}\sum_{k_1=1}^{\lfloor\omega_1\lfloor v_1T/2\rfloor\rfloor}\sum_{\substack{k_2=1\\\eqref{res2}}}^{\lfloor\omega_2\lfloor v_2T/2\rfloor\rfloor}\sum_{l,m,n,o=0}^{\infty}|o-l|\psi_l\psi_m\psi_n\psi_o\frac{1}{2k_1\lfloor v_2T/2\rfloor-2k_2\lfloor v_1T/2\rfloor}\nonumber\\
\leq&\frac{Cv_1v_2}{T}\sum_{k_1=1}^{\lfloor\omega_1\lfloor v_1T/2\rfloor\rfloor}\sum_{\substack{k_2=1\\\eqref{res2}}}^{\lfloor\omega_2\lfloor v_2T/2\rfloor\rfloor}\frac{1}{2k_1\lfloor v_2T/2\rfloor-2k_2\lfloor v_1T/2\rfloor},
\end{align*}
which by simple calculations is seen to be of order $O(\log(T)/T)=O(1/T^{1-\epsilon'})$. The same arguments show that we can change the condition in the summation over $p_2$ to $0\leq p_2\leq2\lfloor v_1T/2\rfloor-1$ by comiting an error  of order $O(1/T^{1-\epsilon'})$ and we thus have that \eqref{cov1} is equal to 
\begin{align}
&\frac{v_1v_2}{T}\sum_{k_1=1}^{\lfloor\omega_1\lfloor v_1T/2\rfloor\rfloor}\sum_{k_2=1}^{\lfloor\omega_2\lfloor v_2T/2\rfloor\rfloor}f(\lambda_{k_1,2\lfloor\frac{v_1T}{2}\rfloor})f(\lambda_{k_2,2\lfloor\frac{v_2T}{2}\rfloor})\nonumber\\
\times&\frac{1}{4\lfloor\frac{v_1T}{2}\rfloor\lfloor\frac{v_2T}{2}\rfloor}\sum_{p_1,p_2=0}^{2\lfloor\frac{v_1T}{2}\rfloor-1}\exp(-i(\lambda_{k_1,2\lfloor\frac{v_1T}{2}\rfloor}-\lambda_{k_2,2\lfloor\frac{v_2T}{2}\rfloor})(p_1-p_2))+O(\frac{1}{T^{1-\epsilon'}})\label{cov2}
\end{align}
[note that we assumed $v_1 \leq v_2$]. In the next step we intend to replace $\lambda_{k_2,2\lfloor\frac{v_2T}{2}\rfloor}$ in the function $f$ by $\lambda_{k_1,2\lfloor\frac{v_1T}{2}\rfloor}$. By employing a Taylor expansion it is straightforward to see that the resulting error can be up to a constant bounded by
\begin{align}
\frac{v_1v_2}{T}&\sum_{k_1=1}^{\lfloor v_1T/2\rfloor}\sum_{k_2=1}^{\lfloor  v_2T/2\rfloor}\frac{|\lambda_{k_1,2\lfloor\frac{v_1T}{2}\rfloor}-\lambda_{k_2,2\lfloor\frac{v_2T}{2}\rfloor}|}{4\lfloor\frac{v_1T}{2}\rfloor\lfloor\frac{v_2T}{2}\rfloor}\Big|\sum_{p_1,p_2=0}^{2\lfloor\frac{v_1T}{2}\rfloor-1}\exp(-i(\lambda_{k_1,2\lfloor\frac{v_1T}{2}\rfloor}-\lambda_{k_2,2\lfloor\frac{v_2T}{2}\rfloor})(p_1-p_2))\Big|\nonumber\\
=&O(\frac{1}{T^{1-\epsilon'}}),\label{approx1}
\end{align}
where the proof of the last equality will be given at a later stage. This yields that, in  \eqref{cov2}, we can replace $\lambda_{k_2,2\lfloor\frac{v_2T}{2}\rfloor}$ in the argument of the spectral density $f$ with $\lambda_{k_1,2\lfloor\frac{v_1T}{2}\rfloor}$  by making an error of order $O(1/T^{1-\epsilon'})$, and we next want to show that the sum over all $(k_1,k_2)$ with $|k_12\lfloor\frac{v_2T}{2}\rfloor-k_22\lfloor\frac{v_1T}{2}\rfloor|\geq v_1v_2a_T$ for some increasing sequence $a_T$ is negligible in the limit. For this reason, we first note that, by employing the geometric series formula and the identity $|1-\exp(ix)|^2=4\sin^2(x/2)$, we get
\begin{align}
&\Big|\sum_{p_1,p_2=0}^{2\lfloor\frac{v_1T}{2}\rfloor-1}\exp(-i(\lambda_{k_1,2\lfloor\frac{v_1T}{2}\rfloor}-\lambda_{k_2,2\lfloor\frac{v_2T}{2}\rfloor})(p_1-p_2))\Big| \nonumber\\
=&\frac{\Big|\big(1-\exp(-i(\lambda_{k_1,2\lfloor\frac{v_1T}{2}\rfloor}-\lambda_{k_2,2\lfloor\frac{v_2T}{2}\rfloor})2\lfloor\frac{v_1T}{2}\rfloor)\big)\big(1-\exp(i(\lambda_{k_1,2\lfloor\frac{v_1T}{2}\rfloor}-\lambda_{k_2,2\lfloor\frac{v_2T}{2}\rfloor})2\lfloor\frac{v_1T}{2}\rfloor)\big)\Big|}{|1-\exp(-i(\lambda_{k_1,2\lfloor\frac{v_1T}{2}\rfloor}-\lambda_{k_2,2\lfloor\frac{v_2T}{2}\rfloor}))|^2} \nonumber\\
\leq& \sin^{-2}\Big((\lambda_{k_1,2\lfloor\frac{v_1T}{2}\rfloor}-\lambda_{k_2,2\lfloor\frac{v_2T}{2}\rfloor})/2\Big). \label{geometricseries}
\end{align}

If we now use this bound in combination with $\sin(\pi x)\geq Cx$ for $x\in[0,1/2]$, assume that $a_T/T^2 \rightarrow 0$ and $ T^{1+\alpha}/a_T \rightarrow 0$ hold for some $0<\alpha <1/2$, and write $(*)_\geq^-$ for the condition $|k_12\lfloor\frac{v_2T}{2}\rfloor-k_22\lfloor\frac{v_1T}{2}\rfloor|\geq v_1v_2a_T$, we obtain
\begin{align*}
&\frac{v_1v_2}{T}\sum_{k_1=1}^{\lfloor\omega_1\lfloor \frac{v_1T}{2}\rfloor\rfloor}f^2(\lambda_{k_1,2\lfloor\frac{v_1T}{2}\rfloor})\frac{1}{4\lfloor\frac{v_1T}{2}\rfloor\lfloor\frac{v_2T}{2}\rfloor}\sum_{\substack{k_2=1\\ (*)_\geq^- }}^{\lfloor\omega_2\lfloor v_2T/2\rfloor\rfloor}\Big|\sum_{p_1,p_2=0}^{2\lfloor\frac{v_1T}{2}\rfloor-1}\exp(-i(\lambda_{k_1,2\lfloor\frac{v_1T}{2}\rfloor}-\lambda_{k_2,2\lfloor\frac{v_2T}{2}\rfloor})(p_1-p_2))\Big|\\
\leq&\frac{v_1v_2C}{T}\sum_{k_1=1}^{\lfloor\omega_1\lfloor \frac{v_1T}{2}\rfloor\rfloor}4\lfloor\frac{v_1T}{2}\rfloor\lfloor\frac{v_2T}{2}\rfloor\sum_{\substack{k_2=1\\ (*)_\geq^-}}^{\lfloor\omega_2\lfloor v_2T/2\rfloor\rfloor}
\frac{1}{(k_12\lfloor\frac{v_2T}{2}\rfloor-k_22\lfloor\frac{v_1T}{2}\rfloor)^2}\\
\leq&\frac{v_1v_2C}{T}\sum_{k_1=1}^{\lfloor\omega_1\lfloor \frac{v_1T}{2}\rfloor\rfloor}\frac{\lfloor v_2T/2\rfloor}{\lfloor v_1T/2\rfloor}\sum_{\substack{k_2=1\\ (*)_\geq^-}}^{\lfloor\omega_2\lfloor v_2T/2\rfloor\rfloor}
\frac{1}{(k_1\frac{\lfloor v_2T/2\rfloor}{\lfloor v_1T/2\rfloor}-k_2)^2}\leq\frac{Cv_2^2}{T}\sum_{k_1=1}^{\lfloor\omega_1\lfloor \frac{v_1T}{2}\rfloor\rfloor}\sum_{k_2=v_2T^{\alpha}}^{\infty}\frac{1}{k_2^2}=O(\frac{1}{T^{\alpha(1-\epsilon')}})
\end{align*}
for some positive constant $C$ and any $\epsilon'>0$. This implies that $V_{T,1}$ is the same as
\begin{align}
&\frac{1}{T}\sum_{k_1=1}^{\lfloor\omega_1\lfloor \frac{v_1T}{2}\rfloor\rfloor}f^2(\lambda_{k_1,2\lfloor\frac{v_1T}{2}\rfloor})A_T(v_1,v_2,k_1,\omega_2) +O(\frac{1}{T^{\alpha(1-\epsilon')}})\label{vorcov4} 
\end{align}
with
\begin{align*}
A_T(v_1,v_2,k_1,\omega_2):=\frac{v_1v_2}{4\lfloor\frac{v_1T}{2}\rfloor\lfloor\frac{v_2T}{2}\rfloor}\sum_{\substack{k_2=1\\ (*)_\leq^-}}^{\lfloor\omega_2\lfloor v_2T/2\rfloor\rfloor}\sum_{p_1,p_2=0}^{2\lfloor\frac{v_1T}{2}\rfloor-1}\exp(-i\frac{2\pi (k_12\lfloor\frac{v_2T}{2}\rfloor-k_22\lfloor\frac{v_1T}{2}\rfloor)}{4\lfloor\frac{v_1T}{2}\rfloor \lfloor\frac{v_2T}{2}\rfloor}(p_1-p_2)) \nonumber
\end{align*}
and where $(*)_\leq^-$ stands for the condition $|k_12\lfloor\frac{v_2T}{2}\rfloor-k_22\lfloor\frac{v_1T}{2}\rfloor|\leq v_1v_2 a_T$.
If we now consider some additional sequence $b_T$  satisfying $T^{1/2+\alpha}/b_T\rightarrow 0$ and $b_T/T^{1-\alpha}\rightarrow 0$ for the same $\alpha>0$ as in the definition of the sequence $a_T$, we get that \eqref{vorcov4} is the same as
\begin{align}
\frac{1}{T}\sum_{k_1=b_T}^{\lfloor\omega_1\lfloor \frac{v_1T}{2}\rfloor\rfloor-b_T}f^2(\lambda_{k_1,2\lfloor\frac{v_1T}{2}\rfloor})A_T(v_1,v_2,k_1,\omega_2)+O(\frac{1}{T^{\alpha(1-\epsilon')}})\label{cov4}.
\end{align}
The property 
\begin{align*}
A_T(v_1,v_2,k_1,\omega_2)=
\begin{cases}
v_1v_2+O(\frac{1}{T^{\alpha}})\quad &\text{if }k_1\in\{b_T,..., \lfloor\omega_2\lfloor \frac{v_1T}{2}\rfloor\rfloor-b_T\}\\
O(\frac{1}{T^{\alpha}})\quad &\text{if }k_1\in\{\lfloor\omega_2\lfloor \frac{v_1T}{2}\rfloor\rfloor+b_T,...,\lfloor \frac{v_1T}{2}\rfloor\}
\end{cases},
\end{align*}
which is obtained by tedious but straightforward calculations [where we basically use the geometric formula, $|1-\exp(ix)|^2=4\sin^2(x/2)$ and a Taylor expansion in combination with the growth rates on $a_T$ and $b_T$], now implies 
with  $v_1\leq v_2$ that
\begin{align}
V_{1,T}=&\frac{\min(v_1,v_2)v_1 v_2}{2\pi}\int_0^{\min(\omega_1,\omega_2)\pi}f^2(\lambda)d\lambda +O(\frac{1}{T^{\alpha(1-\epsilon')}})
\label{varianzfast}
\end{align}
for any $\epsilon'>0$. So, concerning the term $V_{T,1}$, it remains to prove the assertion \eqref{approx1}, and for this purpose we consider the division
\begin{align*}
&\frac{v_1v_2}{T}\sum_{k_1=1}^{\lfloor v_1T/2\rfloor}\sum_{k_2=1}^{\lfloor  v_2T/2\rfloor}\frac{|\lambda_{k_1,2\lfloor\frac{v_1T}{2}\rfloor}-\lambda_{k_2,2\lfloor\frac{v_2T}{2}\rfloor}|}{4\lfloor\frac{v_1T}{2}\rfloor\lfloor\frac{v_2T}{2}\rfloor}\Big|\sum_{p_1,p_2=0}^{2\lfloor\frac{v_1T}{2}\rfloor-1}\exp(-i(\lambda_{k_1,2\lfloor\frac{v_1T}{2}\rfloor}-\lambda_{k_2,2\lfloor\frac{v_2T}{2}\rfloor})(p_1-p_2))\Big|\\
=&S_{\eqref{res1}}+S_{\eqref{res2}}
\end{align*}
into two separate terms, where $S_{\eqref{res1}}$ and $S_{\eqref{res2}}$ denote the sums over all $(k_1,k_2)$ satisfying conditions \eqref{res1} and \eqref{res2} respectively. We now show that each of these sums vanishes as $T\rightarrow\infty$. For $S_{\eqref{res1}}$ we observe that, for fixed $k_1$, there exist at most two values for $k_2$ such that $\eqref{res1}$ holds and the respective summand is non-vanishing. For each of these values it furthermore holds $|\lambda_{k_1,2\lfloor\frac{v_1T}{2}\rfloor}-\lambda_{k_2,2\lfloor\frac{v_2T}{2}\rfloor}|\leq 2\pi (2\lfloor v_2T/2\rfloor)^{-1}$ and we thus obtain
$$S_{\eqref{res1}}\leq\frac{v_1v_2}{T}\sum_{k_1=1}^{\lfloor v_1T/2\rfloor}\frac{2 \pi}{\lfloor v_2T/2\rfloor}=O(\frac{1}{T}).$$
By using  \eqref{geometricseries}  we see that $S_{\eqref{res2}}$ is not larger than
\begin{align}
\frac{v_1v_2}{T}\sum_{k_1=1}^{\lfloor v_1T/2\rfloor}\sum_{\substack{k_2=1\\\eqref{res2}}}^{\lfloor  v_2T/2\rfloor}\frac{|2k_1\lfloor v_2T/2\rfloor-2k_2\lfloor v_1T/2\rfloor|}{(4\lfloor v_1T/2\rfloor\lfloor v_2T/2\rfloor)^2}\sin^{-2}(\frac{\pi(2k_1\lfloor v_2T/2\rfloor-2k_2\lfloor v_1T/2\rfloor)}{4\lfloor v_1T/2\rfloor\lfloor v_2T/2\rfloor}) 
\end{align}
and because of  $\sin{\pi x}\geq Cx$ for $x\leq1/2$    this term is  bounded by 
\begin{align*}
&\frac{v_1v_2}{T}\sum_{k_1=1}^{\lfloor v_1T/2\rfloor}\sum_{\substack{k_2=1\\\eqref{res2}}}^{\lfloor  v_2T/2\rfloor} |2k_1\lfloor v_2T/2\rfloor-2k_2\lfloor v_1T/2\rfloor|^{-1}=\frac{v_1v_2}{2T\lfloor v_1T/2\rfloor}\sum_{k_1=1}^{\lfloor v_1T/2\rfloor}\sum_{\substack{k_2=1\\\eqref{res2}}}^{\lfloor  v_2T/2\rfloor} |k_1\frac{\lfloor v_2T/2\rfloor}{\lfloor v_1T/2\rfloor}-k_2|^{-1} .
\end{align*}
If we  use the fact that, because of \eqref{res2}, $|k_1\frac{\lfloor v_2T/2\rfloor}{\lfloor v_1T/2\rfloor}-k_2|\geq 1$  holds and that it is  $\sum_{k=1}^n k^{-1} \sim\log(n)$, we obtain that this term is of order $O(v_1v_2\log(\lfloor v_2T/2\rfloor)T^{-1})=O(\log(T)/T)$, implying \eqref{approx1}.

Similar arguments can be used to show that the second term $V_{T,2}$ [arising from the application of the product theorem in \eqref{cov1}] vanishes as $T\rightarrow\infty$. More precisely, the condition $\text{cum}(Z_i,Z_j)=\delta_{ij}$ in this case implies that $p_1=p_2+l-n$ and $q_1=q_2-o+m$ have to hold, which by the same arguments as in the calculation of $V_{T,1}$ yields that $V_{T,2}$ is up to an error term of order $O(T^{-\alpha})$ bounded by
\begin{align}
&\frac{Cv_1v_2}{T}\sum_{k_1=1}^{\lfloor\omega_1\lfloor \frac{v_1T}{2}\rfloor\rfloor} \frac{1}{4\lfloor\frac{v_1T}{2}\rfloor\lfloor\frac{v_2T}{2}\rfloor}\sum_{\substack{k_2=1\\ (*)_{\not=}^+}}^{\lfloor\omega_2\lfloor v_2T/2\rfloor\rfloor}\sum_{p_2,q_2=0}^{2\lfloor\frac{v_1T}{2}\rfloor-1} \sin^{-2}(\frac{\pi (k_12\lfloor\frac{v_2T}{2}\rfloor+k_22\lfloor\frac{v_1T}{2}\rfloor)}{4\lfloor\frac{v_1T}{2}\rfloor \lfloor\frac{v_2T}{2}\rfloor})
\label{VT2.1}
\end{align}
with $(*)_{\not=}^+$ standing for the condition $|k_12\lfloor \frac{v_2T}{2}\rfloor+k_22\lfloor \frac{v_1T}{2}\rfloor| \neq4\lfloor \frac{v_1T}{2}\rfloor\lfloor \frac{v_2T}{2}\rfloor$ being satisfied. By splitting the sum over $k_1$ and $k_2$ up into  summands satisfying 
\begin{align}
&|k_12\lfloor v_2T/2\rfloor+k_22\lfloor v_1T/2\rfloor|\leq 2\lfloor v_1T/2\rfloor\lfloor v_2T/2\rfloor\label{condi1}\\
2\lfloor v_1T/2\rfloor\lfloor v_2T/2\rfloor<&|k_12\lfloor v_2T/2\rfloor+k_22\lfloor v_1T/2\rfloor|\leq4\lfloor v_1T/2\rfloor\lfloor v_2T/2\rfloor-2\lfloor v_1T/2\rfloor \label{condi2}
\end{align}
and employing  the bounds $\sin(\pi x)\geq Cx$ for $x\in[0,1/2]$ and $\sin(\pi x)\geq 2(1-x)$ for $x\in[1/2,1]$ it can be shown by similar arguments as provided in the calculation of the term $V_{T,1}$ that \eqref{VT2.1} is of order $O(\log(T)/T)$. The stated covariance kernel for the process $\{\sqrt{T}\hat D_{T}(v,\omega)\}_{(v,\omega)\in[0,1]^2}$ is then obtained with \eqref{varianzfast} using linearity arguments. \\

\textit{Proof of (iii):} For the treatment of the cumulants of order $l\geq3$ we employ the representation
\be
\hat D_T(v,\omega)=\frac{1}{T}\sum_{j=1}^T\sum_{k=1}^{j/2}\phi_{v,\omega,T}(j,\lambda_{k,j})I_j(\lambda_{k,j}),
\ee
where the functions $\phi_{v,\omega,T}$ are defined by
\be
\phi_{v,\omega,T}(j,\lambda):=v1_{\{2\lfloor vT/2\rfloor\}}(j)1_{[0,\frac{2\pi \lfloor\omega\lfloor vT/2\rfloor\rfloor}{2\lfloor vT/2\rfloor}]}(\lambda)-v^21_{\{2\lfloor T/2\rfloor\}}(j)1_{[0,\frac{2\pi \lfloor\omega\lfloor T/2\rfloor\rfloor}{2\lfloor T/2\rfloor}]}(\lambda)
\ee
for $v, \omega \in [0,1]$. By applying the product theorem for cumulants [cf. Theorem 2.3.2 in \cite{brillinger1981}], the Gaussianity and independence of the innovation process $\{Z_t\}_{t\in\Z}$ and the definition of the functions $\phi_{v,\omega,T}$   analogously to the proof of Theorem 6.1 d) in \cite{preussC}, we obtain
\begin{align*}
\text{cum}(\sqrt{T}\hat D_T(v_1,\omega_1),...,\sqrt{T}\hat D_T(v_l,\omega_l))=\sum_{\nu}V(\nu),
\end{align*}
where 
\begin{align*}
V(\nu):=&\frac{1}{T^{l/2}}\sum_{(j_1,...,j_l)\in A_{T}(v_1,...,v_l)}\sum_{k_1=1}^{j_1/2}\hdots \sum_{k_l=1}^{j_l/2}\sum_{m_1,...,m_l=0}^{\infty}\sum_{n_1,...,n_l=0}^{\infty}\prod_{s=1}^l\frac{1}{j_s}\sum_{p_1,q_1=0}^{j_1-1}\hdots\sum_{p_l,q_l=0}^{j_l-1}\\
&\times \prod_{s=1}^l[\phi_{v_s,\omega_s,T}(j_s,\lambda_{k_s,j_s})\psi_{m_s}\psi_{n_s}\exp(-i\lambda_{k_s,j_s}(p_s-q_s))]\text{cum}(Y_{a,b};(a,b)\in\nu_1)\hdots\text{cum}(Y_{a,b};(a,b)\in\nu_l),
\end{align*}
with $Y_{i,1}:=Z_{1+p_i-m_i}$, $Y_{i,2}:=Z_{1+q_i-n_i}$,  $A_{T}(v_1,...,v_l):=\{2\lfloor\frac{v_1T}{2}\rfloor,2\lfloor\frac{T}{2}\rfloor\}\times\hdots\times\{2\lfloor\frac{v_lT}{2}\rfloor,2\lfloor\frac{T}{2}\rfloor\}$
and where the summation is performed over all indecomposable partitions $(\nu_1,...,\nu_l)$  of the scheme
\begin{equation}\label{schema}
\begin{matrix}
Y_{1,1}&Y_{1,2}\\
\vdots&\vdots\\
Y_{l,1}&Y_{l,2}\\
\end{matrix}
\end{equation}
having exactly two elements in each set $\nu_i$. Without loss of generality we restrict ourselves to the partition
$\bar{\nu}:=\cup_{i=1}^{l-1}(Y_{i,1},Y_{i+1,2})\cup(Y_{l,1},Y_{1,2})$
and by simple calculations we obtain
\begin{align}
V(\bar\nu)=&\frac{1}{T^{l/2}}\sum_{(j_1,...,j_l)\in A_{T}(v_1,...,v_l)}\sum_{k_1=1}^{j_1/2}\hdots \sum_{k_l=1}^{j_l/2}\prod_{s=1}^l\phi_{v_s,\omega_s,T}(j_s,\lambda_{k_s,j_s})\sum_{m_1,...,m_l=0}^{\infty}\sum_{n_1,...,n_l=0}^{\infty}\prod_{s=1}^l\frac{1}{j_s}\sum_{\substack{p_1=0\\\eqref{cond2}}}^{j_1-1}\sum_{\substack{p_2=0\\\eqref{cond2}}}^{j_2-1}\hdots\sum_{\substack{p_l=0\\\eqref{cond1}}}^{j_l-1}\nonumber\\
&\times \prod_{s=1}^l[\psi_{m_s}\psi_{n_s}]\times\exp(-i\lambda_{k_1,j_1}(p_1-p_l+m_l-n_1))\times\prod_{s=2}^l\exp(-i\lambda_{k_s,j_s}(p_s-p_{s-1}+m_{s-1}-n_s))\label{V1},
\end{align}
where the conditions 
\begin{align}
0\leq &q_1=p_l-m_l+n_1\leq j_1 \label{cond1}\\
0\leq &q_{i+1}=p_i-m_i+n_{i+1}\leq j_{i+1}\quad i=1,...,l-1\label{cond2}
\end{align}
follow from the independence of the $Z_i$ and the specific form of the partition $\bar\nu$. \eqref{cond1} and \eqref{cond2} imply that $|n_1-m_l|\leq \max(j_1,j_l)$ and  $|n_{i+1}-m_i|\leq \max(j_i,j_{i+1})$ for $i=1,...,l-1$ have to hold. By additionally applying the bound
\begin{align*}
\Big|\frac{1}{j_s}\sum_{k_s=1}^{j_s/2}\phi_{v_s,\omega_s,T}(j_s,\lambda_{k_s,j_s})\exp(-i\lambda_{k_s,j_s}(r))\Big|\leq\frac{C}{r \text{ mod } j_s/2},
\end{align*} 
which holds uniformly in $j_s$, $v,\omega$ and $r$ mod $j_s/2\neq 0$ [see (A.2) in \cite{eichler2008}] it can be seen that \eqref{V1} is bounded by
\begin{align}
&\frac{C^l}{T^{l/2}}\sum_{\substack{(j_1,...,j_l)\in \\A_{T}(v_1,...,v_l)}}\sum_{m_1,...,m_l=1}^{\infty}\sum_{\substack{n_1,...,n_l=1\\|n_1-m_l|\leq\max(j_1,j_l)\\|n_{i+1}-m_i|\leq\max(j_i,j_{i+1})}}^{\infty}\prod_{s=1}^l|\psi_{m_s}\psi_{n_s}|\sum_{p_1=0}^{j_1-1}\sum_{\substack{p_2=0\\|p_2-p_1+m_1-n_2|<\frac{j_2}{2}}}^{j_2-1}\sum_{\substack{p_3=0\\|p_3-p_2+m_2-n_3|<\frac{j_3}{2}}}^{j_3-1}\nonumber\\
\times&\hdots\sum_{\substack{p_l=0\\|p_1-p_l+m_l-n_1|<\frac{j_1}{2}\\|p_l-p_{l-1}+m_{l-1}-n_l|<\frac{j_l}{2}}}^{j_l-1}\prod_{s=1}^l\frac{1}{|p_s-p_{s-1}+m_{s-1}-n_s|}\times\prod_{s=1}^l1(p_s\notin\{z_{s1},z_{s2}\})\label{V2},
\end{align}
where we identify 0 with $l$ and $l+1$ with 1 and define $z_{s1}:=p_{s-1}-m_{s-1}+n_s$ and $z_{s2}:=p_{s+1}+m_s-n_{s+1}$ and exploit the fact that the sums over $|p_s-p_{s-1}+m_{s-1}-n_s|\geq j_s/2$ and $p_s\in\{z_{s1},z_{s2}\}$ are of smaller or the same order. We now define the sets $A_i:=[0,j_i-1]\setminus\{[z_{i1}-1,z_{i1}+1]\cup [z_{i2}-1,z_{i2}+1]\}$, $i=1,...,l$ and obtain that \eqref{V2} is bounded through
\begin{align*}
\frac{C^l}{T^{l/2}}&\sum_{\substack{(j_1,...,j_l)\in \\A_{T}(v_1,...,v_l)}}\sum_{m_1,...,m_l=1}^{\infty}\sum_{\substack{n_1,...,n_l=1\\|n_1-m_l|\leq\max(j_1,j_l)\\|n_{i+1}-m_i|\leq\max(j_i,j_{i+1})}}^{\infty}\prod_{s=1}^l|\psi_{m_s}\psi_{n_s}|\\
\times&\sum_{p_1=0}^{j_1-1}\int_{A_2\times\hdots\times A_l}\prod_{s=1}^l\frac{1}{|p_s-p_{s-1}+m_{s-1}-n_s|}\times\prod_{s=1}^l1(p_s\notin\{z_{s1},z_{s2}\})d(p_2,...,p_l),
\end{align*}
which, by following the proof of Theorem 6.1 d) in \cite{preussC} for  $D=0$ and exploiting  $\text{card}(A_{T}(v_1,...,v_l))=2^l$, is bounded by
\begin{align*}
\frac{C^l}{T^{l/2}}&\sum_{\substack{(j_1,...,j_l)\in \\A_{T}(v_1,...,v_l)}}\sum_{m_1,...,m_l=1}^{\infty}\sum_{\substack{n_1,...,n_l=1\\|n_1-m_l|\leq T\\|n_{i+1}-m_i|\leq T\\m_1-n_1+...+m_l-n_l\neq 0}}^{\infty}\prod_{s=1}^l|\psi_{m_s}\psi_{n_s}|\frac{1}{|m_1-n_1+m_2-n_2+\hdots+m_l-n_l|}&\sum_{p_1=0}^{j_1-1}\log(2lT)^{l-1}\\
=&C^lO(\frac{\log(T)^l}{T^{l/2-1}}).
\end{align*}
Part (iii) then follows due to the fact that the number of partitions $\nu$ of the table \eqref{schema} is finite.

\textbf{Proof of \eqref{Gleichstetigkeit}:} Because of the decomposition
\begin{align*}
\hat D_T(v,\omega)=\frac{v}{T}\sum_{k=1}^{\lfloor \omega   \lfloor \frac{vT}{2}\rfloor   \rfloor} I_{2\lfloor \frac{vT}{2}\rfloor}\big(\lambda_{k,\lfloor\frac{ vT}{2}\rfloor} \big)-\frac{v^2}{T}\sum_{k=1}^{\lfloor \omega  \lfloor\frac{T}{2}\rfloor   \rfloor} I_{2\lfloor\frac{T}{2}\rfloor}\big(\lambda_{k,2\lfloor\frac{T}{2}\rfloor}\big)=:\hat D_{T,1}(v,\omega)-\hat D_{T,2}(v,\omega),
\end{align*}
it is sufficient to show asymptotic stochastic equicontinuity for the processes $\{\sqrt{T}\hat{ D}_{T,1}(v,\omega)\}_{(v,\omega)\in[0,1]^2}$ and $\{\sqrt{T}\hat{ D}_{T,2}(v,\omega)\}_{(v,\omega)\in[0,1]^2}$ seperately. For the sake of brevitiy we restrict ourselves to a proof of \eqref{Gleichstetigkeit} for the first summand and note that asymptotic stochastic equicontinuity for the second summand can be obtained by employing the fact $\hat D_{T,2}(v,\omega)=v D_{T,1}(1,\omega)$. Therefore we use the representation 
\be\label{representation}
\hat D_{T,1}(v,\omega)=\frac{1}{T}\sum_{j=1}^T\sum_{k=1}^{\lfloor  j/2\rfloor}\phi^{(1)}_{v,\omega,T}(j,\lambda_{k,j})I_{j}(\lambda_{k,j}),
\ee
where the functions $\phi^{(1)}_{v,\omega,T}$ are defined by $\phi^{(1)}_{v,\omega,T}(j,\lambda_{k,j}):=v1_{\{2\lfloor vT/2\rfloor\}}(j)1_{[0,\frac{2\pi \lfloor\omega\lfloor vT/2\rfloor\rfloor}{2\lfloor vT/2\rfloor}]}(\lambda_{k,j})$. We furthermore define $d_{T}((v_1,\omega_1),(v_2,\omega_2)):=(|v_1-v_2|+|\omega_1-\omega_2|)^{\beta/2}$ which is a semimetric on the set 
$$\mathcal{P}_T:=\Big\{(v,\omega)\Big|v\in \{\frac{2}{T},\frac{4}{T},...,\frac{2\lfloor T /2\rfloor}{T}\},\omega\in\{0,\frac{1}{\lfloor vT/2\rfloor},...,1-\frac{1}{\lfloor vT/2\rfloor},1\}\Big\}.$$
This notation implies that the equality
\begin{align}
&P\Big( \sup_{\substack{y_1,y_2 \in [0,1]^2:\\ d(y_1,y_2) <\delta }} \sqrt{T}|\hat{ D}_{T,1}(y_1)-\hat{  D}_{T,1}(y_2)|> \eta \Big) = P\Big( \sup_{\substack{y_1,y_2 \in \mathcal{P}_T:\\ d_T(y_1,y_2) <\delta}} \sqrt{T}|\hat{ D}_{T,1}(y_1)-\hat{  D}_{T,1}(y_2)|> \eta \Big) \label{boundP}
\end{align}
holds. In order to show \eqref{Gleichstetigkeit} for the process $\{\sqrt{T}\hat{ D}_{T,1}(v,\omega)\}_{(v,\omega)\in[0,1]^2}$ it is therefore sufficient to show  \eqref{boundP} becomes arbitrarily small for a sufficiently small $\delta$ and large values of $T$. To achieve this goal we denote by $C(u,d_T,\mathcal{P}_T)$ the covering numbers of the set $\mathcal{P}_T$ with respect to $d_T$ and consider the covering integral $J_T(\kappa):=\int_0^{\kappa}[\log\Big(\frac{48C(u,d_T,\mathcal{P}_T)^2}{u}\Big)]^2 du$ [see \cite{wellnervandervaart} for a definition and more details about covering numbers]. We  now proceed by showing the following two assertions:
\begin{itemize}
\item[(i)] For the covering integral it holds
\be\label{condi}
\lim\limits_{\kappa \rightarrow 0}\lim\limits_{T \rightarrow \infty}J_T(\kappa)=0.
\ee
\item[(ii)] For the increments of the process $\{\hat{ D}_{T,1}(v,\omega)\}_{(v,\omega)\in[0,1]^2}$ it holds
\be\label{condii}
\E\big(T^{k/2}(\hat{ D}_{T,1}(v_1,\omega_1)-\hat{  D}_{T,1}(v_2,\omega_2))^k\big)\leq (2k)!C^kd_T((v_1,\omega_1),(v_2,\omega_2))^k
\ee 
for some constant $C>0$, all $(v_i,\omega_i)\in\mathcal{P}_T$ and all even integers $k\in\N$.
\end{itemize}
A similar string of arguments as provided in Theorem 2.4 of \cite{dahlhaus1988} then shows that (i) and (ii) imply that \eqref{boundP} becomes arbitrarily small for $T$ sufficiently large [see the proof of Theorem 2.1 in \cite{detprevet2011b} for more details in a similar framework].\\

\textit{Proof of (i):} The definition of $d_T$ implies that for the covering numbers $C(u,d_T,\mathcal{P}_T)$ the bound \linebreak{$C(u,d_T,\mathcal{P}_T)\leq  \frac{C}{u^{4/\beta}}$} holds, which implies
\begin{align*}
J_T(\kappa)=&\int_0^{\kappa}\Big[\log\Big(48C(u,d_{T},\mathcal{P}_T)^2 u^{-1}\Big)\Big]^2 du\leq C\int_0^{\kappa}[\log(48 u^{-8/\beta-1})]^2 du\\
=&C\Big[\int_0^{\kappa}\log(48)^2du-2\int_0^{\kappa}\log(48)\log(u^{8/\beta+1})du+\int_0^{\kappa}(\log(u^{8/\beta+1}))^2du\Big]\rightarrow 0
\end{align*}
as $\kappa\rightarrow 0$

\textit{Proof of (ii):} The arguments provided in Theorem 2.4 of \cite{dahlhaus1988} show that \eqref{condii} is a consequence of the inequality
\be\label{boundcum}
\text{cum}_l\big(\sqrt{T}(\hat{ D}_{T,1}(v_1,\omega_1)-\hat{  D}_{T,1}(v_2,\omega_2))\big)\leq (2l)!C^ld_T((v_1,\omega_1),(v_2,\omega_2))^l,
\ee
for all $(v_1,\omega_1),(v_2,\omega_2)\in\mathcal{P}_T$ and $l\in\N$.
In order to show \eqref{boundcum} we consider the cases $l=2$ and $l\geq 3$ separately [note that the case $l=1$ is treated by using the same argumentation]. We take $(v_1,\omega_1)\neq (v_2,\omega_2)\in\mathcal{P}_T$, assume without loss of generality that $v_1\leq v_2$, $\omega_1\leq\omega_2$ and define $\phi(j,\lambda):=\phi^{(1)}_{v_1,\omega_1,T}(j,\lambda)-\phi^{(1)}_{v_2,\omega_2,T}(j,\lambda)$. For the case $l=2$ we obtain by following the calculation of the covariances [with $\alpha=\beta$] and applying the linearity of the covariance that $\text{cum}_2(\sqrt{T}(\hat{ D}_{T,1}(v_1,\omega_1)-\hat{  D}_{T,1}(v_2,\omega_2)))$  is equal to
\begin{align*}
&\frac{1}{2\pi}\Big(\int_0^{\omega_1\pi}f^2(\lambda)d\lambda(v_1^3+v_2^3-2v_1^2v_2)+v_2^3\int_{\omega_1\pi}^{\omega_2\pi}f^2(\lambda)d\lambda\Big)+O(\frac{1}{T^{\beta}})\\
\leq &\frac{1}{2\pi}\Big(\int_0^{\omega_1\pi}f^2(\lambda)d\lambda v_2(v_2^2-v_1^2)+\int_{\omega_1\pi}^{\omega_2\pi}f^2(\lambda)d\lambda\Big)+O(\frac{1}{T^{\beta}})\\
\leq&C(|v_1-v_2|+|\omega_1-\omega_2|)+(|v_1-v_2|+|\omega_1-\omega_2|)^{\beta}\leq C(|v_1-v_2|+|\omega_1-\omega_2|)^{\beta},
\end{align*}
where we used that for $(v_1,\omega_1)\neq (v_2,\omega_2)\in\mathcal{P}_T$ the inequality $T^{-1}\leq C(|v_1-v_2|+|\omega_1-\omega_2|)$ holds and that there exists a  constant $C$ such that the bound $x\leq C_1x^{\beta}$ holds uniformly in $x\in[0,2]$.
For $l\geq 3$ we obtain, by following the steps in the   calculation of the higher oder cumulants  in the proof of \eqref{fidis}, and the fact that the number of partitions of the table \eqref{schema} is bounded by $(2l)!2^l$ that 
\begin{align*}
\text{cum}_l\big(\sqrt{T}(\hat{ D}_{T,1}(v_1,\omega_1)-\hat{  D}_{T,1}(v_2,\omega_2))\big)=(2l)!C^lO(\frac{\log(T)^{l-1}}{T^{l/2-1}})=(2l)!C^lO(\frac{1}{T^{l(1/2-1/3-\kappa)}})
\end{align*}
uniformly in $(v_1,\omega_1),(v_2,\omega_2)\in\mathcal{P}_T$ and for any $\kappa>0$. By choosing $\kappa=1/6-\beta/2$  and again considering $T^{-1} \leq C(|v_1-v_2|+|\omega_1-\omega_2|)$ for $(v_1,\omega_1)\neq (v_2,\omega_2)\in\mathcal{P}_T$ we thus obtain
\begin{align*}
\text{cum}_l\big(\sqrt{T}(\hat{ D}_{T,1}(v_1,\omega_1)-\hat{  D}_{T,1}(v_2,\omega_2))\big) \leq (2l)!C^l d_T((v_1,\omega_1),(v_2,\omega_2))^l,
\end{align*}
which completes the proof of Theorem \ref{hauptsatz}. \\

{\bf Proof of part b):} We obtain completely analogously to the proof of part a) [with a little bit more   notation] that $\sup_{v,\omega}\sqrt{T}|\hat D_T(v,\omega)-\E(D_T(v,\omega))|=O_P(1)$. Therefore the claim follows with \eqref{expectation1}. $\hfill \Box$

%BEWEIS BOOTSTRAP
\subsection{Proof of Theorem \ref{theorembootstrap}}
For this proof we again restrict ourselves to the case $d=1$. It follows from \eqref{ARparameterRate} and Lemma 2.3 in \cite{kreisspappol2011} that, for sufficiently large $T$, there exists a MA$(\infty)$ representation
\begin{equation}\label{gefittetARinMA}
X_{t,T}^*=\sum_{l=0}^{\infty}\hat{\psi}_l^{AR}(p)Z_{t-l}^*
\end{equation}
for the sequence $\{X_{t,T}^*\}_{t=1,...,T}$ of bootstrap replicates. By noting that the bootstrapped version \linebreak{$\{\sqrt{T}\hat D_T^*(v,\omega)\}_{(v,\omega)\in[0,1]^2}$} is defined in the same way as the original empirical process \linebreak{$\{\sqrt{T}\hat D_T(v,\omega)\}_{(v,\omega)\in[0,1]^2}$}, with the original time series data $\{X_{t,T}\}_{t=1,...,T}$ being replaced by $\{X_{t,T}^*\}_{t=1,...,T}$, and  that \eqref{gefittetARinMA} corresponds to a stationary MA$(\infty)$ time series model under the nullhypothesis as well as under the alternative, we can perform the proof along the arguments as were provided in the proof of Theorem \ref{hauptsatz}. Therefore we substitute the coefficents $\{\hat{\psi}_l^{AR}(p)\}_{l\geq 0}$, which are now random, for the deterministic sequence $\{\psi_l\}_{l\geq 0}$ and substitute $\{Z_t^*\}_{t\in\Z}$ for $\{Z_t\}_{t\in\Z}$. The assertion  follows by showing that the random error terms of the form
\bea
\frac{\big(\sum_{l=0}^{\infty}|\hat{\psi}_l^{AR}(p)|\big)^{q_1}\big(\sum_{m=0}^{\infty}|m||\hat{\psi}_m^{AR}(p)|\big)^{q_2}}{T}
\eea
$q_1,q_2\in\N$, which take the place of the deterministic errors of the form 
\bea
\frac{\big(\sum_{l=0}^{\infty}|\psi_l|\big)^{q_1}\big(\sum_{m=0}^{\infty}|m||\psi_m(p)|\big)^{q_2}}{T},
\eea
are of order $O_P(1/T)$. This however is implied by the fact $\big(\sum_{l=0}^{\infty}|\hat{\psi}_l^{AR}(p)|\big)^{q_1}\big(\sum_{m=0}^{\infty}|m||\hat{\psi}_m^{AR}(p)|\big)^{q_2}=O_P(1)$ for all $q_1,q_2\in\N$, which is shown in the proof of Theorem 3.2 in \cite{detprevet2011b}. $\hfill \Box$

\bibliographystyle{apalike}

\begin{thebibliography}{}

\bibitem[Akaike, 1973]{akaike1973}
Akaike, H. (1973).
\newblock {\em Information theory and an extension of the maximum likelihood
  principle}.
\newblock Budapest, Akademia Kiado, 267-281.

\bibitem[Anderson, 1971]{anderson1971}
Anderson, T.~W. (1971).
\newblock {\em The Statistical Analysis of Time Series}.
\newblock John Wiley and Sons, New York.

\bibitem[Berg et~al., 2010]{bergpappolitis2010}
Berg, A., Paparoditis, E., and Politis, D.~N. (2010).
\newblock A bootstrap test for time series linearity.
\newblock {\em Journal of Statistical Planning and Inference}, 140:3841--3857.

\bibitem[Brillinger, 1981]{brillinger1981}
Brillinger, D.~R. (1981).
\newblock {\em Time Series: Data Analysis and Theory}.
\newblock McGraw Hill, New York.

\bibitem[Brockwell and Davis, 1991]{brodav1991}
Brockwell, P.~J. and Davis, R.~A. (1991).
\newblock {\em Time Series: Theory and Methods}.
\newblock Springer Verlag, New York.

\bibitem[Chandler and Polonik, 2006]{chandler2006}
Chandler, G. and Polonik, W. (2006).
\newblock Discrimination of locally stationary time series based on the excess
  mass functional.
\newblock {\em Journal of the American Statistical Association}, 101:240û253.

\bibitem[Changli et~al., 2009]{changli2009}
Changli, H., Terõsvirta, T., and Gonzßlez, A. (2009).
\newblock Testing parameter constancy in stationary vector autoregressive
  models against continuous change.
\newblock {\em Econometric Reviews}, 28:225--245.

\bibitem[Chen et~al., 2010]{motivation2}
Chen, Y., Hõrdle, W., and Pigorsch, U. (2010).
\newblock Localized realized volatility modeling.
\newblock {\em Journal of the American Statistical Association},
  105(492):1376--1393.

\bibitem[Dahlhaus, 1988]{dahlhaus1988}
Dahlhaus, R. (1988).
\newblock Empirical spectral processes and their applications to time series
  analysis.
\newblock {\em Stochastic Process and their Applications}, 30:69--83.

\bibitem[Dahlhaus, 1997]{dahlhaus1997}
Dahlhaus, R. (1997).
\newblock Fitting time series models to nonstationary processes.
\newblock {\em Annals of Statistics}, 25(1):1--37.

\bibitem[Dahlhaus, 2000]{dahlhaus2000}
Dahlhaus, R. (2000).
\newblock A likelihood approximation for locally stationary processes.
\newblock {\em Annals of Statistics}, 28(6):1762--1794.

\bibitem[Dahlhaus, 2009]{dahlhaus2009}
Dahlhaus, R. (2009).
\newblock Local inference for locally stationary time series based on the
  empirical spectral measure.
\newblock {\em Journal of Econometrics}, 151:101--112.

\bibitem[Dahlhaus, 2012]{dahlhaus2011}
Dahlhaus, R. (2012).
\newblock Locally stationary processes.
\newblock {\em Handbook of Statistics}, 30.

\bibitem[Dahlhaus and Polonik, 2009]{dahlpolo2009}
Dahlhaus, R. and Polonik, W. (2009).
\newblock Empirical spectral processes for locally stationary time series.
\newblock {\em Bernoulli}, 15:1--39.

\bibitem[Dette et~al., 2011]{detprevet2010}
Dette, H., Preu{\ss}, P., and Vetter, M. (2011).
\newblock A measure of stationarity in locally stationary processes with
  applications to testing.
\newblock {\em Journal of the American Statistical Association},
  106(495):1113--1124.

\bibitem[Dwivedi and Subba~Rao, 2011]{rao2010}
Dwivedi, Y. and Subba~Rao, S. (2011).
\newblock A test for second order stationarity of a time series based on the
  discrete fourier transform.
\newblock {\em Journal of Time Series Analysis.}, 32(1):68--91.

\bibitem[Eichler, 2008]{eichler2008}
Eichler, M. (2008).
\newblock Testing nonparametric and semiparametric hypotheses in vector
  stationary processes.
\newblock {\em Journal of Multivariate Analysis}, 99:968--1009.

\bibitem[Fan, 1994]{thresholding2}
Fan, J. (1994).
\newblock Test of significance based on wavelet thresholding and neyman's
  truncation.
\newblock {\em Journal of the American Statistical Association},
  91(434):674--688.

\bibitem[Fryzlewicz and Ombao, 2009]{piotr2009}
Fryzlewicz, P. and Ombao, H. (2009).
\newblock Consistent classification of nonstationary time series using
  stochastic wavelet representations.
\newblock {\em Journal of the American Statistical Association}, 104:299--312.

\bibitem[Huang et~al., 2004]{huang2004}
Huang, H.-Y., Ombao, H., and Stoffer, D. (2004).
\newblock Discrimination and classification of nonstationary time series using
  the slex model.
\newblock {\em Journal of the American Statistical Association}, 99:763û774.

\bibitem[Jentsch and Subba~Rao, 2012]{rao2012}
Jentsch, C. and Subba~Rao, S. (2012).
\newblock A test for second order stationarity of a multivariate time series.
\newblock Technical report.

\bibitem[Krei{\ss}, 1988]{Kreiss1988}
Krei{\ss}, J.-P. (1988).
\newblock {\em Asymptotic statistical inference for a class of stochastic
  processes}.
\newblock Habilitations\-schrift, Fachbereich Mathematik, Universit{\"a}t
  Hamburg.

\bibitem[Krei{\ss} and Paparoditis, 2012]{krepap2012}
Krei{\ss}, J.~P. and Paparoditis, E. (2012).
\newblock The hybrid wild bootstrap for time series.
\newblock {\em Journal of the American Statistical Association},
  107(499):1073--1084.

\bibitem[Krei{\ss} et~al., 2012]{kreisspappol2011}
Krei{\ss}, J.-P., Paparoditis, E., and Politis, D.~N. (2012).
\newblock On the range of the validity of the autoregressive sieve bootstrap.
\newblock {\em Annals of Statistics}, 39(4):2103.

\bibitem[Neumann and von Sachs, 1997]{neumsach1997}
Neumann, M.~H. and von Sachs, R. (1997).
\newblock Wavelet thresholding in anisotropic function classes and applications
  to adaptive estimation of evolutionary spectra.
\newblock {\em Annals of Statistics}, 25:38--76.

\bibitem[Paparoditis, 2009]{paparoditis2009}
Paparoditis, E. (2009).
\newblock Testing temporal constancy of the spectral structure of a time
  series.
\newblock {\em Bernoulli}, 15:1190--1221.

\bibitem[Paparoditis, 2010]{paparoditis2010}
Paparoditis, E. (2010).
\newblock Validating stationarity assumptions in time series analysis by
  rolling local periodograms.
\newblock {\em Journal of the American Statistical Association},
  105(490):839--851.

\bibitem[Paparoditis and Preu{\ss}, 2013]{pappreuss}
Paparoditis, S. and Preu{\ss}, P. (2013).
\newblock On local properties of frequency domain based tests for stationarity.
\newblock Technical report.

\bibitem[Preu{\ss} and Vetter, 2013]{preussC}
Preu{\ss}, P. and Vetter, M. (2013).
\newblock On discriminating between long-range dependence and non stationarity.
\newblock {\em Electronic Journal of Statistics}, 7:2241--2297.

\bibitem[Preu{\ss} et~al., 2012]{detprevet2011b}
Preu{\ss}, P., Vetter, M., and Dette, H. (2012).
\newblock A test for stationarity based on empirical processes.
\newblock {\em to appear in Bernoulli}.

\bibitem[Preu{\ss} et~al., 2013]{preusssemi}
Preu{\ss}, P., Vetter, M., and Dette, H. (2013).
\newblock Testing semiparametric hypotheses in locally stationary processes.
\newblock {\em Scandinavian Journal of Statistics}, 40.

\bibitem[Priestley and Subba~Rao, 1969]{priestleystat}
Priestley, M. and Subba~Rao, T. (1969).
\newblock A test for non-stationarity of time series.
\newblock {\em Journal of the Royal Statistical Society, Series B},
  31:140--149.

\bibitem[Sakiyama and Taniguchi, 2003]{Taniguchi}
Sakiyama, K. and Taniguchi, M. (2003).
\newblock Testing composite hypotheses for locally stationary processes.
\newblock {\em Journal of Time Series Analysis}, 24(4):483--504.

\bibitem[Sakiyama and Taniguchi, 2004]{sakitani2004}
Sakiyama, K. and Taniguchi, M. (2004).
\newblock Discriminant analysis for locally stationary processes.
\newblock {\em Journal of Multivariate Analysis}, 90:282--300.

\bibitem[Sergides and Paparoditis, 2009]{sergpapa2009}
Sergides, M. and Paparoditis, E. (2009).
\newblock Frequency domain tests of semiparametric hypotheses for locally
  stationary processes.
\newblock {\em Scandinavian Journal of Statistics}, 36:800--821.

\bibitem[Starica and Granger, 2005]{starica2005}
Starica, C. and Granger, C. (2005).
\newblock Nonstationarities in stock returns.
\newblock {\em The Review of Economics and Statistics}, 87:503--522.

\bibitem[van~der Vaart and Wellner, 1996]{wellnervandervaart}
van~der Vaart, A. and Wellner, J. (1996).
\newblock {\em Weak Convergence and Empirical Processes}.
\newblock Springer, Berlin.

\bibitem[von Sachs and Neumann, 2000]{vonSachs2000}
von Sachs, M. and Neumann, M.~H. (2000).
\newblock A wavelet-based test for stationarity.
\newblock {\em Journal of time series analysis}, 21:597--613.

\bibitem[Whittle, 1951]{whittle1}
Whittle, P. (1951).
\newblock {\em Hypothesis Testing in Time Series Analysis}.
\newblock HUppsala: Almqvist and Wiksell.

\end{thebibliography}

\end{document}